\documentclass[journal]{IEEEtran}

\IEEEoverridecommandlockouts                              
\let\labelindent\relax
\usepackage{amsthm, amssymb, amsmath, enumitem, pdfpages, color, cite, caption, graphicx, relsize, subfigure}
\newtheorem{thm}{Theorem}		    		 
\newtheorem{cor}{Corollary}		    		 
\newtheorem{definition}{Definition}	    		 
\newtheorem{lem}{Lemma}		    		 
	     	    		 
\newtheorem{rmk}{Remark}	               		 
\newtheorem{assn}{Assumption}	    		 
\def\mc{\mathcal}
\def\mb{\mathbf}
\def\mbb{\mathbb}
\def\ra{\rightarrow}
\def\bs{\boldsymbol}

\title{Optimization over time-varying directed graphs with row and column-stochastic matrices}
\author{Fakhteh~Saadatniaki,~\IEEEmembership{Student~Member,~IEEE}, Ran~Xin,~\IEEEmembership{Student~Member,~IEEE},\\and~Usman~A.~Khan$^\dagger$,~\IEEEmembership{Senior~Member,~IEEE}%
\thanks{$^\dagger$The authors are with the Department of Electrical and Computer Engineering at Tufts University,~Medford,~MA~02155; {\tt\footnotesize \{fakhteh.saadatniaki@,ran.xin@,khan@ece.\}tufts.edu}. This work is partially supported by an NSF Career award: CCF~\#~1350264.}}

\begin{document}
\maketitle

\begin{abstract}
In this paper, we provide a distributed optimization algorithm, termed as~TV-$\mc{AB}$,  that minimizes a sum of convex functions over time-varying, random directed graphs. Contrary to the existing work, the algorithm we propose does not require eigenvector estimation to estimate the (non-$\mb{1}$) Perron eigenvector of a stochastic matrix. Instead, the proposed approach relies on a novel information mixing approach that exploits both row- and column-stochastic weights to achieve agreement towards the optimal solution when the underlying graph is directed. We show that~TV-$\mc{AB}$ converges linearly to the optimal solution when the global objective is smooth and strongly-convex, and the underlying time-varying graphs exhibit bounded connectivity, i.e., a union of every~$C$ consecutive graphs is strongly-connected. We derive the convergence results based on the stability analysis of a linear system of inequalities along with a matrix perturbation argument. Simulations confirm the findings in this paper. 

\IEEEkeywords Optimization, distributed algorithms, first-order methods, time-varying graphs, directed graphs.

\end{abstract}

\section{Introduction}
With the advent of~5G, promising higher bandwidth and faster data rates, emerging technologies like the Internet of Things, self-driving cars, and smart devices are coming to the forefront. In these applications, it is of paramount interest to learn hidden parameters from the data collected by the individual units~\cite{6515173}. For example, self-driving cars rely heavily on computer vision in order to correctly identify pedestrians, highway lanes, or traffic signs in all light and weather conditions. Problems such as these can be framed as classification, regression, or risk minimization, at the heart of which is a simple sum of cost optimization. The sheer size of data and privacy concerns limit data sharing and thus solutions to the underlying optimization problems must be developed that are local and distributed~\cite{NIPS2003_2448,NIPS2010_4164,NIPS2013_4976}. However, departing from the traditional approaches, what is needed now are algorithms that are applicable to mobile and autonomous agents resulting in a time-varying and non-deterministic information exchange. 

Recently, there has been a large body of work on distributed optimization, where the goal is to minimize a sum of costs:
\begin{align*}
\min_{\mb{x}} \frac{1}{n}\textstyle \sum_{i=1}^n f_i(\mb{x}), 
\end{align*}
such that each objective function,~$f_i:\mbb{R}^p\ra\mbb{R}$, is private to agent~$i$. In order to solve this problem, the agents exchange information with nearby nodes over a sparse communication graph. Such problems arise for example in sensor networks~\cite{distributed_Rabbit}, large-scale machine learning~\cite{7219480, wai2018multi}, distributed estimation~\cite{exactdiffusion1}, and localization~\cite{safavi2018distributed}. When the graphs are static and undirected, early work on first-order methods include~\cite{DOPT1, uc_Nedic, duchi2012dual} with the convergence rate of~${O(\frac{\ln k}{\sqrt{k}})}$ for arbitrary convex functions and~${O(\frac{\ln k}{k})}$ for strongly-convex functions, where~$k$ denotes the number of iterations. The sublinear convergence is due to the use of diminishing step-sizes. The rate improves to linear with a constant step-size but at the expense of a sub-optimal solution~\cite{DGD_Yuan}. Methods based on Lagrangian dual~\cite{ADMM_Wei, ADMM_Mota} converge faster but suffer from a high computational burden as they require solving a subproblem at each iteration. 

%\vspace{-0.025cm}
Optimization over directed graphs is developed in~\cite{opdirect_Tsianous, opdirect_Tsianous2, opdirect_Nedic}, where push-sum~\cite{ac_directed0, ac_directed} is used to achieve consensus among the agents. The convergence rate, with diminishing step-sizes, is~${O(\frac{\ln k}{\sqrt{k}})}$ for arbitrary convex functions and~${O(\frac{\ln k}{k})}$ for strongly-convex functions. In contrast, Refs.~\cite{D-DGD,D-DPS} use an alternate approach called surplus consensus~\cite{ac_Cai1} to achieve consensus but with the same convergence rates as~\cite{opdirect_Tsianous, opdirect_Tsianous2, opdirect_Nedic}. The main reason for slow convergence is that a local gradient is used at each agent, which requires diminishing step-sizes to ensure optimality. To overcome this challenge, Refs.~\cite{AugDGM,harness, diging} replace the local gradient with an estimate of the global gradient, with the help of dynamic consensus~\cite{DAC} over undirected graphs, and show linear convergence to the optimal solution. This gradient estimation approach was combined with push-sum (type) methods in~\cite{add-opt, diging, linear_row,FROST} to achieve linear convergence to the optimal solution over directed graphs. Related work along these lines also include~\cite{EXTRA} over undirected graphs and~\cite{DEXTRA} over directed graphs that are related to~\cite{add-opt, diging}. An alternate approach that builds on~\cite{D-DGD,D-DPS} and does not use push-sum has been recently developed in~\cite{AB}, where a row- and a column-stochastic matrix are used simultaneously to achieve linear convergence to the optimal over directed graphs. Accelerated methods can be found in~\cite{dnesterov, nic_mom, ABm}, while non-convex problems are considered in~\cite{ asy_sonata}.

%\vspace{-0.05cm}
In this paper, our focus is on time-varying directed graphs. Early work on time-varying communication among the agents can be found in~\cite{cc_nedic, cc_Lobel2}, where the exchange is undirected and in~\cite{opdirect_Nedic3}, where the exchange is directed. These methods are built on local gradients at each agent and are thus sublinear. Recent work includes~\cite{lu2018geometrical}, where the authors establish a geometrically converging distributed optimization algorithm over directed graphs under uncoordinated yet bounded step-sizes, and~\cite{xu2018convergence}, where agents communicate over graphs subject to random link failures. Work on random graphs can be found in~\cite{cc_lobel}, where the problem of constrained convex optimization is investigated for non-differentiable costs under Markovian communication model. For random networks modeled by a sequence of independent, identically distributed~(IID) random matrices drawn from the set
of symmetric, stochastic matrices with positive diagonals,~\cite{6665045} proposes two accelerated distributed Nesterov-like gradient methods featuring resiliency to link failures, reduced computational load, and improved convergence rates compared to other gradient methods. Ref.~\cite{sgd_random} establishes a convergence rate of~$O(\frac{1}{k})$ for distributed stochastic gradient methods over temporally IID random undirected graphs for strongly-convex costs when local gradients are subject to noise that is IID in time and has a finite second moment. Furthermore, asynchronous multi-agent optimization is considered in~\cite{7849202}, where the authors adapt the curvature estimation technique of the Broyden-Fletcher-Goldfarb-Shanno quasi-Newton optimization method~\cite{1976Sgcp,doi:10.1137/0724077, NocedalJorge2006No} for use in asynchronous distributed settings over undirected graphs~\cite{7903733}. An asynchronous implementation of subgradient-push~\cite{opdirect_Nedic} algorithm is also recently developed in~\cite{assran2018asynchronous}.

In this paper, we use gradient estimation as was used in~\cite{AugDGM, harness, add-opt, diging} and extend the~$\mc{AB}$ algorithm introduced in~\cite{AB} to time-varying, random directed graphs. Of relevance in this context is Ref.~\cite{diging}, which uses gradient estimation and push-sum consensus~\cite{ac_directed0, ac_directed} to implement distributed optimization over time-varying graphs. However, the push-sum based methods~\cite{add-opt, diging,linear_row, FROST} involve estimating the (non-$\mb{1}$) Perron eigenvector of a stochastic matrix and an appropriate scaling with the estimated eigenvector components. The eigenvector estimation adds conservatism to the overall algorithm as its convergence may dominate the overall rate. In contrast, the~$\mc{AB}$ algorithm~\cite{AB} does not require such eigenvector estimation and employs both row- and column-stochastic weights in a novel way. The time-varying algorithm proposed in this paper, termed as~TV-$\mc{AB}$, thus is applicable to time-varying, directed graphs without the need of eigenvector estimation resulting from push-sum. 

The~TV-$\mc{AB}$ algorithm we introduce involves a unique and a rather counter-intuitive way of mixing information among the agents. As the graph at each iteration may not be strongly-connected, the mechanics of~TV-$\mc{AB}$ can be explained over a single directed edge,~$i\longrightarrow j$. First, we note that~{TV-$\mc{AB}$} involves two state updates: one with a row-stochastic weight matrix,~$\mc{A}$, and the other with a column-stochastic weight matrix,~$\mc B$. The update involving~$\mc{A}$ is standard where the receiving agent~$j$ implements a sum-preserving update to its past and  the incoming information from agent~$i$ , while agent~$i$ assigns a weight of ~$1$ to its past since it does not receive any information. However, the additional update with the column-stochastic weights,~$\mc{B}$, requires the transmitting agent~$i$ to implement a strictly stable update (by assigning a weight less than~$1$ to its past) and the receiving agent to implement an unstable update (by assigning weights that sum to~$>1$ to its past and the incoming information from agent~$i$) in order to maintain column-stochasticity of~$\mc{B}$. In other words, the updates involving~$\mc{B}$ are not sum-preserving unlike the traditional information fusion.

We show that~TV-$\mc{AB}$ converges linearly to the optimal solution when each local objective is smooth and the global objective is strongly-convex. The graph at each iteration can be generated randomly in an arbitrary fashion as long as the union of every~$C$ consecutive graphs is strongly-connected. This notion is known as bounded connectivity and is standard in the consensus and optimization literature on time-varying graphs, see e.g.,~\cite{diging, opdirect_Nedic}. The bounded connectivity notion enables us to obtain more concrete convergence results as without this assumption, the analysis is restricted to the expected behavior of the optimization algorithm, see e.g.,~\cite{cc_Lobel2,6665045,xu2018convergence}. We show linear convergence with the help of a linear system of inequalities along with a matrix perturbation argument. 

We now describe the rest of the paper. Section~\ref{s:pf} formulates the problem and introduces the assumptions necessary to algorithm development. Section~\ref{ss:alg} develops and interprets the time-varying~$\mc{AB}$ algorithm. Details on the convergence analysis are presented in Section~\ref{s:conv_analysis} while Section~\ref{s:lin_conv} states the main result. Finally, Section~\ref{s:sims} provides numerical experiments and Section~\ref{s:conc} contains the concluding remarks.

\textbf{Notation:} We denote column vectors by lowercase bold letters,~$\mb{x}$, and matrices by uppercase italics,~$X$. The~$n\times n$ identity matrix and the~${n\textrm{-dimensional}}$ column vector of all ones are denoted by~$I_n$ and~$\mb{1}_n$, respectively, where the dimension subscripts are dropped if clear from the context. We denote by~$X\otimes Y$, the Kronecker product of two matrices,~$X$ and~$Y$, while~$\rho(X)$ denotes the spectral radius for a matrix,~$X$. For a vector,~$\mb{x}$, its~$i$\textsuperscript{th} element is denoted by~$[\mb{x}]_i$, while for a matrix,~$X$, its~$(i,j)$\textsuperscript{th} element is denoted by~$[X]_{i,j}$. The notation~$\|\cdot\|_2$ denotes the Euclidean norm for vectors and the spectral norm for matrices, whereas~$\|\cdot\|_{\max}$ denotes the~$l_\infty$-norm on the set of square matrices. 

\section{Problem Formulation and Algorithm}\label{s:pf}
In this section, we formulate the distributed optimization problem, state the assumptions, and introduce the~TV-$\mc{AB}$ algorithm. To this aim, we assume that the set of agents,~${\mc{V}=\{1,2,\cdots,n\}}$, communicate according to a time-varying directed graph,~$\mathcal{G}_k(\mathcal{V},\mathcal{E}_k)$, where~$k$ is the discrete-time index and~$\mathcal{E}_k$ is the set of directed communication links at time~$k$. An agent~$j$ can send information to an agent~$i$, i.e.,~${j\longrightarrow i}$, at time~$k$, if and only if~{$(i,j)\in\mc{E}_k$}. The goal of the agents is to collaboratively solve the following problem: 
\begin{align}
\mbox{Problem 1:}\qquad \min_{\mb{x}} f(\mb{x})=\frac{1}{n}\textstyle \sum_{i=1}^n f_i(\mb{x}), 
\end{align}
where each local objective function,~$f_i:\mbb{R}^p\mapsto\mbb{R}$, is held privately at agent~$i$. We next formalize the set of assumptions that are standard in the distributed optimization literature.

\begin{assn}[Strong-convexity]\label{a1:strongly_cvx}
The global objective function~$f$, is~$\mu$-strongly-convex, i.e.,
\begin{align}
f(\mb{x}) \geq f(\mb{y}) + {\nabla}f(\mb{y})^{\top}(\mb{x} - \mb{y}) +\frac{\mu}{2}\|\mb{x}-\mb{y}\|_2^2\label{eq:strongly_cvx}
\end{align}
for any~$\mb{x}$,~$\mb{y}\in\mathbb{R}^p$, where~$\mu > 0$. 
\end{assn}

\noindent For Assumption~\ref{a1:strongly_cvx} to hold it suffices that each~$f_i$ is convex and at least one of them is strongly-convex. Under this assumption, Problem~1 has a unique optimal solution, denoted by~$\mb{x}^{\ast}$.

\begin{assn}[Smoothness]\label{a2:Lipschitz}
Each~$f_i$ is~$\ell_i$-smooth, i.e.,~it is differentiable and has a Lipschitz-continuous gradient. Mathematically, there exists~$\ell_i>0$ such that
\begin{align}
\|{\nabla f}_i(\mb{x})-{\nabla f}_i(\mb{y})\|_2\leq \ell_i\|\mb{x}-\mb{y}\|_2,\label{eq:Lipschitz}
\end{align}
for any~$\mb{x}$,~$\mb{y}\in\mathbb{R}^p$ and~$\forall i\in\mc{V}$.
\end{assn}  
\noindent Assumption~\ref{a2:Lipschitz} implies that~$f=\sum_if_i$ is~$\bar\ell$-smooth, where~${\bar{\ell}=\frac{1}{n}\sum_{i=1}^n\ell_i}$. Furthermore, collecting the local variables in column vectors, i.e.,
\begin{equation*}
\begin{aligned}
\mathbf{x}=&\begin{bmatrix}
{\mb{x}^1} \\
 \vdots \\
 {\mb{x}^n}
\end{bmatrix}, 
& {\mb{f}}(\mathbf{x})=&\begin{bmatrix}
f_1({\mb{x}^1})\\
 \vdots\\
 f_n({\mb{x}^n})
\end{bmatrix},
& {\nabla \mb{f}}(\mathbf{x})=&\begin{bmatrix}
{\nabla }f_1({\mb{x}^1})\\
 \vdots\\
 {\nabla f}_n({\mb{x}^n})
\end{bmatrix},
\end{aligned}
\end{equation*}
we note that~$\mb{f}$ is~$L$-smooth, where~$L=\max_i\{\ell_i\}$.
\begin{assn}[$C$-bounded strong-connectivity%Uniform~$C$-strong-connectivity
]\label{a3:joint_conn}
For the sequence~$\{\mathcal{G}_k=(\mathcal{V},\mathcal{E}_k\subseteq \mathcal{V}\times\mathcal{V})\}$ of time-varying directed graphs, there exists some positive integer~${C}$ such that the aggregate digraph~$\mathcal{G}_k^{{C}}\triangleq (\mathcal{V},\textstyle \cup_{l=k}^{k+C-1}\mathcal{E}_l)$ is strongly-connected~${\forall\,k\geq 0}$.
\end{assn}

\begin{assn}[Weights]\label{a4:weight_mats}
For the sequence~${\{\mathcal{G}_k=(\mathcal{V},\mathcal{E}_k)\}}$ of time-varying directed graphs and the sequences,~$\{A_k\}$ and~$\{B_k\}$, of~$n\times n$ matrices compliant with~$\mathcal{G}_k$, i.e.,~${(i,j)\in \mathcal{E}_k\Leftrightarrow [A_{k}]_{i,j},\,[B_{k}]_{i,j}\neq 0}$, the following hold. 
\begin{enumerate}[label=(\roman*)]
\item Stochasticity:~$\{A_k\}$ and~$\{B_k\}$ are row- and column-stochastic, respectively.
\item Aperiodicity:~$\mathcal{G}_k$ has self-loops; i.e.,~$[A_{k}]_{i,i}>0$ and~$[B_{k}]_{i,i}>0,\,\forall\,i\in\mc{V}$ and ~$\forall\,k\geq 0$.
\item Uniform positivity:~There are scalars~${0<\alpha,\,\beta<1}$ such that~$[A_{k}]_{i,j}\geq\alpha$ and~$[B_{k}]_{i,j}\geq\beta$,~$\forall\,(i,j)\in \mathcal{E}_k, k\geq 0$.
\end{enumerate}
\end{assn}

\noindent The strong-connectivity bound~$C$ introduced in Assumption~\ref{a3:joint_conn} and the uniform positivity bounds~$\alpha$ and~$\beta$ in Assumption~\ref{a4:weight_mats} are not required to be known at any of the agents. They are only used in the analysis of the algorithm. 
% _______________________________________________________  Algorithm Development ________________________________________________________________________ %
\subsection{Algorithm Development}\label{ss:alg}
We now describe the~TV-$\mathcal{AB}$ algorithm to solve Problem~P1. At each time~$k$, agent~$i\in\mc{V}$ maintains two variables,~$\mb{x}_k^{i}$,~$\mb{y}_k^{i}~$, both in~$\mbb{R}^p$, initialized with arbitrary~$\mb{x}_0^{i}~$ and~${\mb{y}_0^{i}={\nabla f}_i(\mb{x}_0^{i})}$. The~$\mb{x}_k^i$-update at each agent is essentially gradient descent, albeit after mixing incoming information, and where the descent direction is given by an estimate of the global gradient,~$\mb{y}_k^i$, instead of the local gradient,~$\nabla f_i(\mb{x}_k^i)$. The~$\mb{y}_k^i$-update at each agent tracks the global gradient and is based off of column-stochastic weights. 

We first use a simple framework to explain the algorithm where only one edge~$i\longrightarrow j$ is active at time~$k$. The~$\mb{x}_k$-update follows the standard sum-preserving notion where weights assigned to the past information are non-negative and sum~to~$1$:
\begin{align*}
\mb{x}_{k + 1}^{i} =& \mb{x}_k^{i}-\eta \mb{y}_k^{i},\qquad \forall\,i\neq j,\\
\mb{x}_{k + 1}^{j} =& [A_k]_{(j,i)}\mb{x}_k^{i} + [A_k]_{(j,j)}\mb{x}_k^{j} - \eta \mb{y}_k^{j},
\end{align*}
where~$\eta$ is a constant step-size. The weight matrix,~$A_k$, behind this update is row-stochastic: Each diagonal element,~$[A_k]_{(i,i)}=1,\forall\, i\neq j$, while the~$j$\textsuperscript{th} row has only two positive elements such that~$[A_k]_{(j,i)}+[A_k]_{(j,j)}=1$. 

Defining the auxiliary variable~${\mb{z}_{k+1}^i}$ as the successive difference,~${{\nabla f}_{{i}}(\mb{x}_{k + 1}^{{i}})-{\nabla f}_{{i}}(\mb{x}_{k}^{\bs{i}})}$, of the gradients, with~${\mb{z}_{0}^i=\mb{0}_p}$, the~$\mb{y}_k$-update is given by
\begin{align*}
\mb{y}_{k + 1}^m =& [B_k]_{(m,m)}\mb{y}_k^m+\mb{z}_{k+1}^m,\qquad m\neq i,m\neq j,\\
\mb{y}_{k + 1}^{{i}} =&[B_k]_{(i,i)}\mb{y}_k^{{i}}+\mb{z}_{k+1}^i,\\
\mb{y}_{k + 1}^{{j}} =&[B_k]_{(j,i)}\mb{y}_k^{{i}}+[B_k]_{(j,j)}\mb{y}_k^{{j}}+\mb{z}_{k+1}^j.
\end{align*}
Note that every column of~$B_k$, except the~$i$\textsuperscript{th}, has only one non-zero element as no agent other~than~$i$ is transmitting.~Hence, for~$B_k$ to be column-stochastic,~${[B_k]_{(j,j)}=1,\forall j\neq i}$. For~the transmitting agent, we have~$[B_k]_{(i,i)}+[B_k]_{(j,i)}=1$. In contrast to the traditional row- or doubly-stochastic updates, the~$\mb{y}_{k+1}^i$-update is locally stable as~$[B_k]_{(i,i)}<1$, while the~$\mb{y}_{k+1}^j$-update is locally unstable as~${[B_k]_{(j,i)}+[B_k]_{(j,j)}>1}$.  

\subsection{The TV-$\mc{AB}$ Algorithm}
The simple scenario discussed above can be generalized to arbitrary graphs, resulting into the following algorithm:
\begin{subequations}\label{eq:alg_agent}
\begin{align}
\mb{x}_{k + 1}^i =& \sum_{j=1}^n [A_{k}]_{i,j}\mb{x}_k^{j}-\eta \mb{y}_k^{i},\label{eq:gd_agent}\\
\mb{y}_{k + 1}^i =& \sum_{j=1}^n [B_{k}]_{i,j}\mb{y}_k^{j}+\mb{z}_k^i,\label{eq:grad_agent}
\end{align}
\end{subequations}
where the weights~$[A_{k}]_{i,j}$ and~$[B_{k}]_{i,j}$ satisfy Assumption~\ref{a4:weight_mats}. Letting~${\mc{A}_k\triangleq A_k\otimes I_p}$ and~$\mc{B}_k\triangleq B_k\otimes I_p$,  we present~Eqs.~\eqref{eq:alg_agent} in a vector-matrix form, where the local variables~$\mb{x}_k^i$’s and~$\mb{y}_k^i$’s and gradients~${\nabla f}_i(\mb{x}_k^i)$ are stored in the column vectors~$\mb{x}_k$,~$\mb{y}_k$, and~${\nabla \mb{f}}(\mb{x}_{k})$, respectively:
\begin{subequations}\label{eq:alg_aggreg} 
\begin{align}
\mathbf{x}_{k+1}&=\mc{A}_k\mathbf{x}_k-\eta \mathbf{y}_k,\label{eq:state}\\
\mathbf{y}_{k+1}&=\mc{B}_k\mathbf{y}_k+\mb{z}_k,  \label{eq:tracking}
\end{align}
\end{subequations} 
where~$\mathbf{z}_k={\nabla \mb{f}}(\mathbf{x}_k)-{\nabla \mb{f}}(\mathbf{x}_{k-1})$, and~${\mathbf{z}}_0=\mb{0}_{np}$. The weight matrices,~$\mc{A}_k$ and~$\mc{B}_k$, are row- and column-stochastic, respectively. However, since each~$\mc{G}_k$ is not necessarily strongly-connected, the weights~$\mc{A}_k$ and~$\mc{B}_k$ are not necessarily primitive or irreducible; thus, the standard Perron-Frobenius arguments are not applicable here. To overcome this isse, we use the notions of absolute probability sequences, ergodicity, and multi-step contractions, a recap of these concepts is provided in the Appendix: Section~\ref{app:prelim}. 

\section{Convergence Analysis}\label{s:conv_analysis}
To proceed with the analysis, we perform a state transformation:~${\mb{s}_k=(V_k^{-1}\otimes I_p)\mb{y}_k}$ on~$\mb{y}_k$, where~${V_k=\textrm{diag}[\mb{v}_k]}$ and~$\mb{v}_k$ follows Eq.~\eqref{eq:left_eig}. The TV-$\mc{AB}$ algorithm is thus equivalently written as 
\begin{subequations}\label{eq:alg_aggreg} 
\begin{align}
\mb{v}_{k+1} &= B_k\mb{v}_k,\label{eq:left_eig}\\
\mathbf{x}_{k+1}&=\mc{A}_k\mathbf{x}_k-\eta (V_k\otimes I_p)\mathbf{s}_k,\label{eq:state}\\
\mathbf{s}_{k+1}&={\mc{R}_k}\mathbf{s}_k+({V}^{-1}_{k+1}\otimes I_p)({\nabla \mb{f}}(\mathbf{x}_{k+1})-{\nabla \mb{f}}(\mathbf{x}_k)),\label{eq:s_tracking}
\end{align}
\end{subequations} 
where~$\mb{v}_0=\mb{1}_n$,~${\mc{R}_k\triangleq R_k\otimes I_p}$, and~${R_k={V}^{-1}_{k+1}{B}_k{V}_k}$. It can be verified that~$\{{R}_k\}$ is a sequence of row-stochastic matrices for which the absolute probability sequence is~${\{\mb{v}_k\}}$; see Appendix: Section~\ref{app:prelim} on absolute probability sequences. 

We now proceed with the convergence analysis of the equivalent algorithm in Eqs.~\eqref{eq:left_eig}-\eqref{eq:s_tracking}. Our approach rests on a few quantities that we describe next: 
\begin{enumerate}[label=(\roman*)]
\item~$\overline{{\mb{x}}}^{\textrm{w}}_k=(\bs{\phi}^{\top}_k\otimes I_p)\mb{x}_k$, which is the average of~$\mb{x}_k^i$'s weighted by the absolute probability sequence,~$\{\bs{\phi}_k\}$, of~$A_k$'s, see Corollary~\ref{cor:phi_cor} in the Appendix: Section~\ref{app:prelim};
\item~$\widetilde{\mathbf{x}}_k^{\textrm{w}}=\mathbf{x}_k-\mb{1}_n\otimes\overline{{\mb{x}}}^{\textrm{w}}_k$, which can be regarded as the weighted consensus error in the network;
\item~$\mathbf{r}_k=\mb{1}_n\otimes\overline{{\mb{x}}}^{\textrm{w}}_k-\mathbf{1}_n\otimes {\mb{x}^{\ast}}$, which is the optimality gap associated with the weighted average;
\item~$\widetilde{\mathbf{s}}^{\textrm{w}}_k={\mathbf{s}}_k-(\mb{1}_n\mb{v}_k^{\top}\otimes I_p){\mathbf{s}}_k$, which is an error term corresponding to gradient estimation.
\end{enumerate}
With the help of these quantities, we define the vector
\begin{align}
\mb{t}_{k}=\begin{bmatrix}
\|\widetilde{\mathbf{x}}^{\textrm{w}}_{k}\|_2\\
\|{\mathbf{r}}_{k}\|_2\\
\|\widetilde{\mathbf{s}}^{\textrm{w}}_{k}\|_2
\end{bmatrix},\label{eq:lin_ineq}
\end{align}
and show that it goes to zero as~$k\ra\infty$. Clearly, if~$\mb{t}_{k}\ra\mb{0}$, then~$\mb{x}_k\ra\mb{1}_n\otimes\mb{x}^*$, and rate of convergence of TV-$\mc{AB}$ is upper bounded by the rate at which~$\mb{t}_{k}\ra\mb{0}$. To establish that~$\mb{t}_{k}\ra\mb{0}$, we derive a linear system of inequalities that expresses the evolution of~$\mb{t}_k$ in the following form:
\begin{align}\label{xyz}
\begin{bmatrix}
\mb{t}_{k+1}\\
\vdots\\
\mb{t}_{k-(\overline{C}-2)}
\end{bmatrix}\leq&{M(\eta)}\begin{bmatrix}
\mb{t}_k\\
\vdots\\
\mb{t}_{k-(\overline{C}-1)}
\end{bmatrix},
\end{align}
where the elements of~$M(\eta)$ are the coefficients of the linear system. Clearly, if~$\rho(M(\eta))<1$, then~$\mb{x}_k\ra\mb{1}_n\otimes\mb{x}^*$ at least at the rate of~$\mc{O}(\rho(M(\eta))^k)$. 

\begin{figure}[!h]
\centering
\includegraphics[width=0.8\columnwidth]{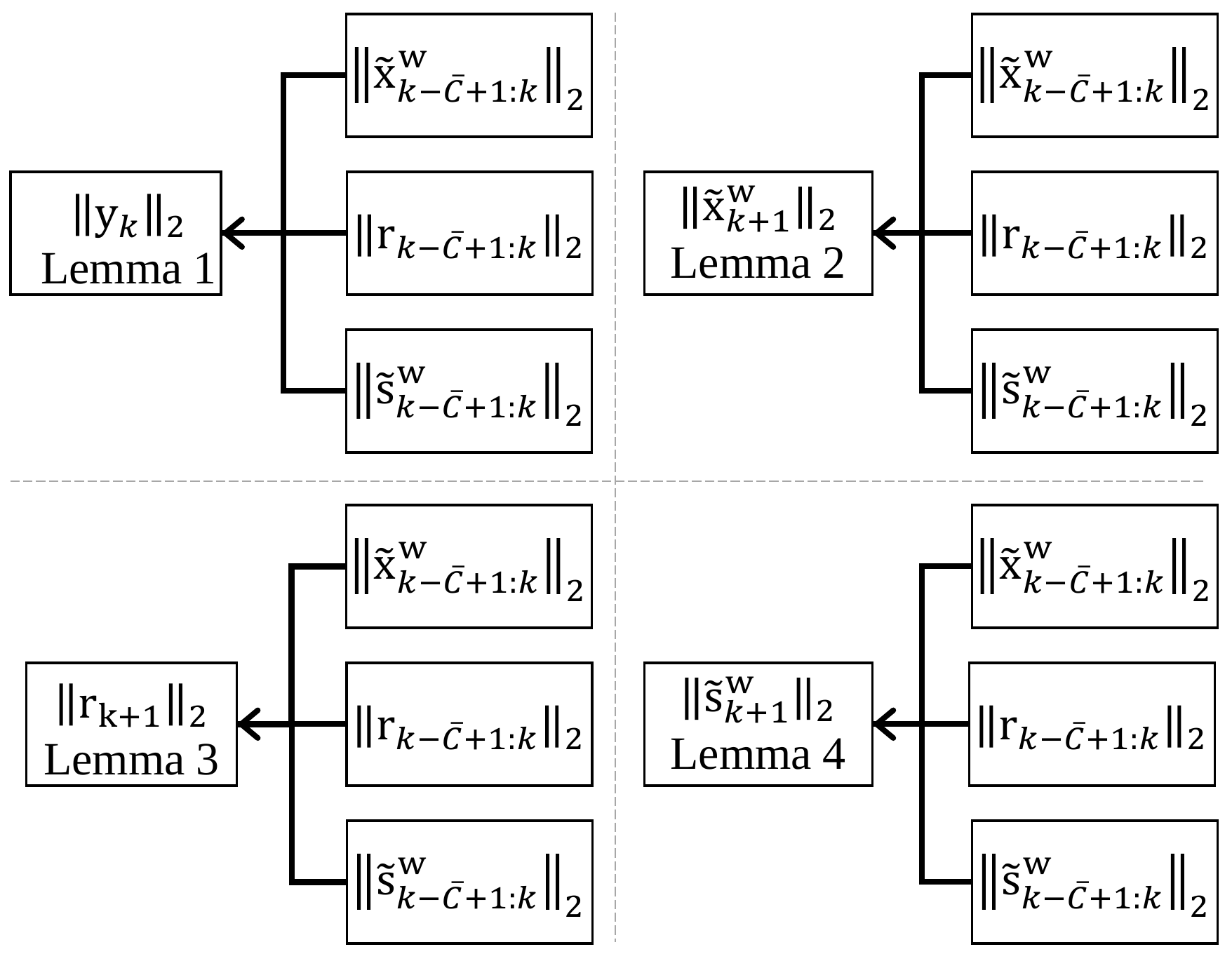}
\caption{Roadmap of deriving the linear system of inequalities.}\label{fig:roadmap}
\end{figure}

Fig.~\ref{fig:roadmap} provides a roadmap to establish Eq.~\eqref{xyz}. The next four lemmas provide the corresponding inequalities, whereas the proofs are deferred to the Appendix. The system of Eq.~\eqref{xyz} is then analyzed in the next subsection.

\begin{rmk}
The constant~${\overline{C}=\max \{\overline{C}_{\mc{A}},\overline{C}_{\mc{B}}\}}$ used in the rest of the paper, ensures concurrent multi-step contractions for both variables~$\widetilde{\mathbf{x}}_k^{\textrm{w}}$ and~$\widetilde{\mathbf{s}}_k^{\textrm{w}}$; see Lemma~\ref{lem:contract_A} and Corollary~\ref{cor:contract_B} in the Appendix: Section~\ref{app:prelim} for more details.
\end{rmk}

\begin{lem} \label{lem:y}The following inequality holds~$\forall\,k\geq 0$:
\begin{align*}
\|\mathbf{y}_k\|_2\leq&{n}L\| \widetilde{\mathbf{x}}^{\textrm{w}}_k\|_2+{n}L\| {\mathbf{r}}_k\|_2+\| \widetilde{\mathbf{s}}^{\textrm{w}}_k\|_2.
\end{align*}
\end{lem}
\begin{IEEEproof}
See Appendix: Section~\ref{app:lem_4}.
\end{IEEEproof}

\begin{lem}\label{lem:x_tilde}The following inequality holds~$\forall\,k\geq \overline{C}-1$:
{\noindent\small \begin{align*}
\begin{aligned}
\|\widetilde{\mathbf{x}}^{\textrm{w}}_{k+1}\|_2\leq &\big({\gamma_{\mc{A}}}+\eta Q_{\mc{A}}{n}L\big)\|\widetilde{\mathbf{x}}^{\textrm{w}}_{k-(\bar{C}-1)}\|_2+\eta Q_{\mc{A}}nL\sum_{l=0}^{\overline{C}-2}\|\widetilde{\mathbf{x}}^{\textrm{w}}_{k-l}\|_2\\
&+\eta Q_{\mc{A}}{n}L\Big(\|{\mathbf{r}}_{k-(\bar{C}-1)}\|_2+\sum_{l=0}^{\overline{C}-2}\|{\mathbf{r}}_{k-l}\|_2\Big)\\
&+\eta Q_{\mc{A}}\Big(\|\widetilde{\mathbf{s}}^{\textrm{w}}_{k-(\bar{C}-1)}\|_2+\sum_{l=0}^{\overline{C}-2}\|\widetilde{\mathbf{s}}^{\textrm{w}}_{k-l}\|_2\Big),
\end{aligned}
\end{align*}}where~$\gamma_{\mc{A}}$ and~$Q_{\mc{A}}$ are the constants defined in Lemma~\ref{lem:contract_A} in the Appendix: Section~\ref{app:prelim}.
\end{lem}
\begin{IEEEproof}
See Appendix: Section~\ref{app:lem_5}.
\end{IEEEproof}

\begin{lem}\label{lem:r}The following inequality holds~$\forall\,k\geq 0$:
\begin{align*}
\begin{aligned}
\|{\mathbf{r}}_{k+1}\|_2
\leq &\eta n L\|\widetilde{\mathbf{x}}^{\textrm{w}}_k\|_2+\Big(1-\eta\frac{{\mu}}{n^{nC-1}}\Big)\|\mb{r}_k\|_2+\eta\sqrt{n} \|\widetilde{\mathbf{s}}^{\textrm{w}}_k\|_2.
\end{aligned}
\end{align*}
\end{lem}
\begin{IEEEproof}See Appendix: Section~\ref{app:lem_6}.
\end{IEEEproof}

\begin{lem}\label{lem:s_tilde}The following inequality holds~$\forall\,k\geq \overline{C}-1$:
{\noindent\small \begin{align*}
\|\widetilde{\mathbf{s}}^{\textrm{w}}_{k+1}\|_2\leq& m\sqrt{n}(2+\eta L)\Big(\|\widetilde{\mb{x}}^{\textrm{w}}_{k-(\overline{C}-1)}\|_2+\sum_{l=0}^{\overline{C}-2}\|\widetilde{\mb{x}}^{\textrm{w}}_{k-l}\|_2\Big)\\
&+\eta m{n}L\Big(\|{\mb{r}}_{k-(\overline{C}-1)}\|_2+\sum_{l=0}^{\overline{C}-2}\|{\mb{r}}_{k-l}\|_2\Big)\\
&+(\eta m+{\gamma_{\mc{B}}})\|\widetilde{\mathbf{s}}^{\textrm{w}}_{k-(\overline{C}-1)}\|_2+\eta m\sum_{l=0}^{\overline{C}-2}\|\widetilde{\mb{s}}^{\textrm{w}}_{k-l}\|_2,
\end{align*}}where~$m=n^{nC}Q_{\mc{B}}L$, and~$\gamma_{\mc{B}}$ and~$Q_{\mc{B}}$ are the constants defined in Corollary~\ref{cor:contract_B} in the Appendix: Section~\ref{app:prelim}.
\end{lem}
\begin{IEEEproof}See Appendix: Section~\ref{app:lem_7}.
\end{IEEEproof}

% Room left for 6 lines!
\subsection{The resulting linear system of inequalities}
Summarizing the results of Lemmas~\ref{lem:x_tilde}-\ref{lem:s_tilde}, for~${0<\eta<\frac{2}{nL}}$, ~Eq.~\eqref{xyz} can be expanded as follows:
{\small\begin{align*}
\mb{t}&_{k+1}\leq\underbrace{\begin{bmatrix}
\eta Q_{\mc{A}}{n}L& \eta Q_{\mc{A}}{n}L& \eta Q_{\mc{A}}\\
\eta nL&1-\eta\frac{{\mu}}{n^{nC-1}}&\eta \sqrt{n}\\
m\sqrt{n}(2+\eta L)&\eta m{n}L&\eta m
\end{bmatrix}}_{M_1}\mb{t}_k\\
&+\underbrace{\begin{bmatrix}
\eta Q_{\mc{A}}nL& \eta Q_{\mc{A}}nL& \eta Q_{\mc{A}}\\
0&0&0\\
m\sqrt{n}(2+\eta L)&\eta m{n}L&\eta m
\end{bmatrix}}_{M_2}(\mb{t}_{k-1}+\ldots+\mb{t}_{k-(\overline{C}-2)})\nonumber\\
&+\underbrace{\begin{bmatrix}
\gamma_{\mc{A}}+\eta Q_{\mc{A}}{n}L& \eta Q_{\mc{A}}{n}L& \eta Q_{\mc{A}}\\
0&0&0\\
m\sqrt{n}(2+\eta L)&\eta m{n}L&\gamma_{\mc{B}}+\eta m
\end{bmatrix}}_{M_{\overline{C}}}\mb{t}_{k-(\overline{C}-1)},
\end{align*}}which is equivalent to
\noindent{\small \begin{align}
\begin{bmatrix}
\mb{t}_{k+1}\\
\mb{t}_k\\
\mb{t}_{k-1}\\
\vdots\\
\mb{t}_{k-(\overline{C}-2)}
\end{bmatrix}\leq\begin{bmatrix}
{M_1}&{M_2}&
\cdots&{M_2}&{M_{\overline{C}}}\\
I\\
&I\\
&&\ddots\\
&&&I
\end{bmatrix}\begin{bmatrix}
\mb{t}_k\\
\mb{t}_{k-1}\\
\vdots\\
\mb{t}_{k-(\overline{C}-2)}\\
\mb{t}_{k-(\overline{C}-1)}
\end{bmatrix}.\label{eq:S}
\end{align}}The system matrix,~$M(\eta)$, in the above can be partitioned as 
{\small\begin{align}
\underbrace{\begin{bmatrix}
M_1^0&\cdots&M^0_2&M_{\overline{C}}^0\\
I\\
&\ddots\\
&&I\\
\end{bmatrix}}_{M^0}
+\eta 
\underbrace{\begin{bmatrix}
M^E_1&\cdots&M^E_2&M_{\overline{C}}^E\\
\mb{0}\\
&\ddots\\
&&\mb{0}\\
\end{bmatrix}}_{M^E},
\label{eq:M_0}
\end{align}}where
{\small \begin{align*}
\begin{aligned}
M_1^0={\begin{bmatrix}
0 & 0 & 0 \\
0&1&0\\
2m\sqrt{n}&0&0
\end{bmatrix}},&\,\,\,M_1^E={\begin{bmatrix}
Q_{\mc{A}}{n}L& Q_{\mc{A}}{n}L& Q_{\mc{A}}\\
nL&-\frac{{\mu}}{n^{nC-1}}&\sqrt{n}\\
m\sqrt{n}L&m{n}L&m
\end{bmatrix}},&\\
M_2^0={\begin{bmatrix}
0&0& 0\\
0&0&0\\
2m\sqrt{n}&0&0
\end{bmatrix}},&\,\,\,M_2^E={\begin{bmatrix}
Q_{\mc{A}}{n}L& Q_{\mc{A}}{n}L& Q_{\mc{A}}\\
 0&0&0\\
m\sqrt{n}L&m{n}L&m
\end{bmatrix}},&\\
M_{\overline{C}}^0={\begin{bmatrix}
\gamma_{\mc{A}}&0&0\\
0&0&0\\
2m\sqrt{n}&0&\gamma_{\mc{B}}
\end{bmatrix}},&\,\,\,M_{\overline{C}}^E={\begin{bmatrix}
Q_{\mc{A}}{n}L& Q_{\mc{A}}{n}L& Q_{\mc{A}}\\
0&0&0\\
m\sqrt{n}L&m{n}L&m
\end{bmatrix}}.&
\end{aligned}
\end{align*}}

Recall that our goal is to establish the geometric decay of~$\mb{t}_k$ in Eq.~\eqref{eq:lin_ineq}. To this purpose, it is sufficient to show~${\rho(M(\eta))<1}$. As a first step, we finish this section with a lemma on the spectral radius of the matrix~$M^0$ in Eq.~\eqref{eq:M_0} and its corresponding eigenvector.
\begin{lem}\label{lem:M_0}
The spectral radius of the matrix~$M^0$ is~$1$ and~${\lambda=1}$ is a simple eigenvalue of~$M^0$. The left and right eigenvectors,~$M_0\mb{u} = \mb{u},\mb{w}^\top M_0 = \mb{w}^\top$, are given by
\begin{align}
\mb{u} &= \mb{1}_{\overline C} \otimes \begin{bmatrix}0&1&0\end{bmatrix}^{\top},\\
\mb{w}^{\top}&=\begin{bmatrix}
0 & 1&0&\cdots&0
\end{bmatrix}.
\end{align}
\end{lem}
\begin{IEEEproof}
See Appendix: Section~\ref{M0_ev}. 
\end{IEEEproof}

\section{Linear Convergence}\label{s:lin_conv}
We now state the main convergence result for TV-$\mc{AB}$.
\begin{thm}
The spectral radius of~$M(\eta)$ is strictly less than~$1$ when~$\eta$ is sufficiently small. Therefore~$\|\mb{x}_k-\mb{1}_n\otimes{\mb{x}^{\ast}}\|_2$ converges to zero  (at least) at the rate of~${O\big(\rho(M(\eta))^k\big)}$.
\end{thm}
\begin{IEEEproof}
From Lemma~\ref{lem:M_0}, let~$q(\eta)$ be the simple eigenvalue of~$M(\eta)$, as a function of~$\eta$, for which~$q(0)=1$. Recall that~$M(\eta)$ can be partitioned as~$M^0+\eta M^E$ from Eq.~\eqref{eq:M_0}. Borrowing a result from matrix perturbation theory~\cite[Theorem~6.3.12]{matrix}, we have that
\[\left.\frac{dq(\eta)}{d\eta}\right|_{\eta=0}=\frac{ \mb{w}^{\top}M^E\mb{u}}{ \mb{w}^{\top}\mb{u}},\]
where~$\mb{u}$ and~$\mb{w}$ are right and left eigenvectors corresponding to the simple eigenvalue,~$q(0)$. From Lemma~\ref{lem:M_0}, it can be verified that~$\mb{w}^{\top}\mb{u}=1$ and~$\mb{w}^{\top}M^E\mb{u}=-{{\mu}}/{n^{nC-1}}<0$, which implies that~$\textstyle \frac{d}{d\eta}q(\eta)$ is negative. Since the eigenvalues are a continuous function of the elements of a matrix, we have that~$q(\eta)$ decreases for a sufficiently small~$\eta$ (slightly increasing from zero) and the theorem follows.
\end{IEEEproof}

\section{Numerical Experiments}\label{s:sims}
This section illustrates the application and performance of the time-varying~$\mc{AB}$ algorithm in a variety of numerical experiments. In the rest of this section, we adopt a simple uniform weighting strategy to construct the row- and column-stochastic weights~$[A_{k}]_{i,j}$ and~$[B_{k}]_{i,j}$:
\begin{align}\label{eq:A_row}
[A_k]_{i,j}=\left\{\begin{array}{rc}
1/{d_{k,{\textrm{in}}}^i},& (i,j)\in \mathcal{E}_k,\\
0,& (i,j)\notin \mathcal{E}_k,
\end{array}\right.
\end{align}
where~${d_{k,{\textrm{in}}}^i}$ is the in-degree of agent~$i$ at time~$k$; and
\begin{align}\label{eq:B_col}
[B_k]_{i,j}=\left\{\begin{array}{rc}
1/{d_{k,\textrm{out}}^j},& (i,j)\in \mathcal{E}_k,\\
0,& (i,j)\notin \mathcal{E}_k,
\end{array}\right.
\end{align}
where~${d_{k,\textrm{out}}^j}$ is the out-degree of agent~$j$ at time~$k$.

\subsection{Distributed binary classification}\label{ex:B}
\begin{figure}[!ht]
\centering
\includegraphics[width=0.76\columnwidth]{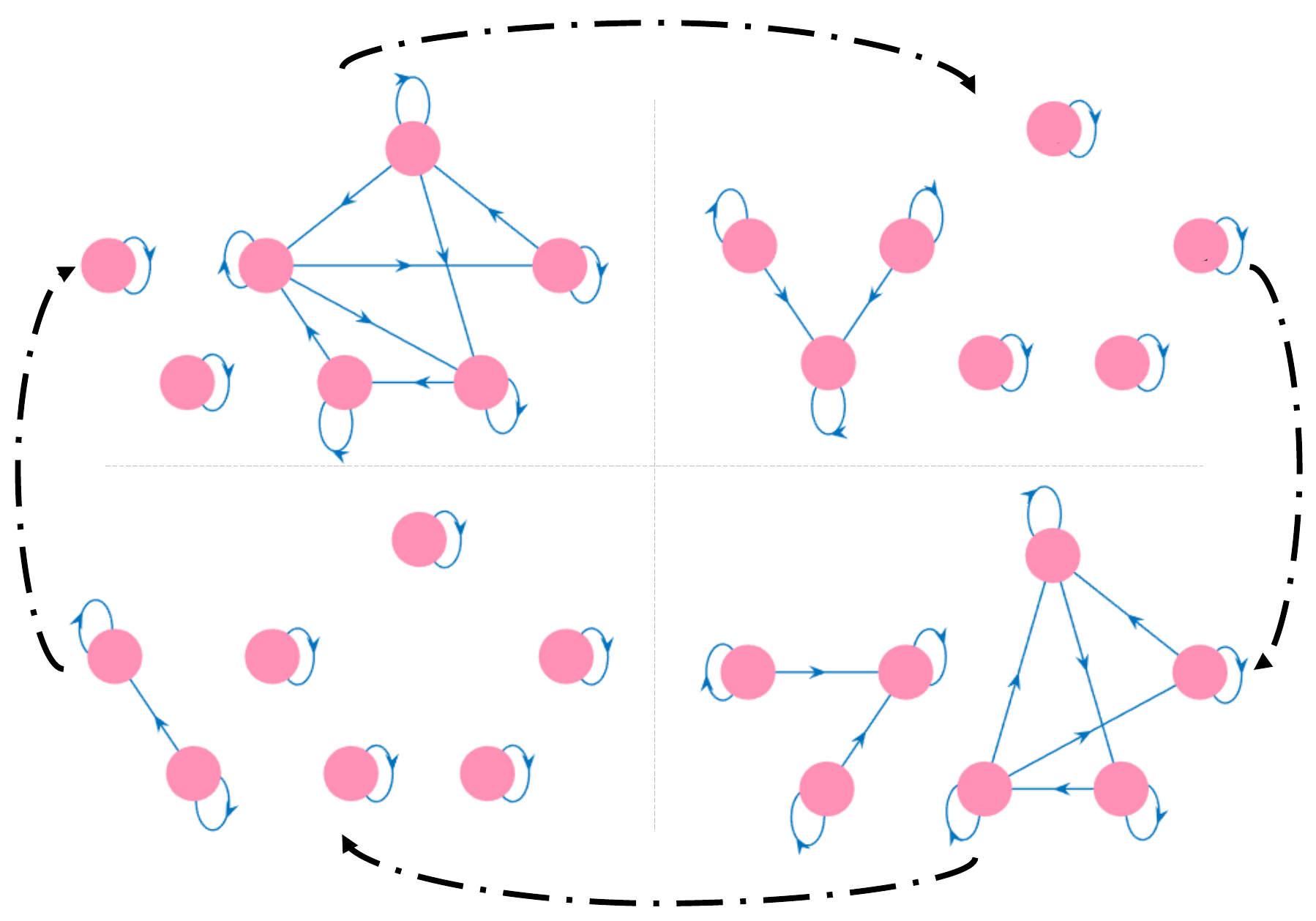}
\caption{Periodic time-varying topology, where the period is~$4$.}
\label{fig:B_1-4}
\end{figure}

In this experiment, we study a binary classification problem using regularized logistic regression. Each agent~$i$ has access to~$m_i$ training samples:~${({\mb{c}^i}^{(j)} ,{y^i}^{(j)})\in\mbb{R}^{p-1}\times \{-1,\,1\}}$, for~${j=1,\,2,\,\cdots,\,m_i}$, where~${\mb{c}^i}^{(j)}$ is the~{$p-1$-dimensional} feature vector of the~$j$\textsuperscript{th} training sample at the~$i$\textsuperscript{th} agent and~${y^i}^{(j)}$ is the corresponding binary label. The agents collaboratively solve the following distributed logistic regression problem:
\begin{align*}
\min_{\mb{w},b}f(\mb{w},b)=&\sum_{i=1}^nf_i(\mb{w},b),
\end{align*} 
where the private loss function~$f_i$ at agent~$i$ is
\begin{align*}
f_i(\mb{w},b)=&\sum_{j=1}^{m_i}\ln \left[1 + e^{\big(-\mb{w}^{\top}{\mb{c}^i}^{(j)}+ b\big){y^i}^{(j)}}\right]+\frac{\lambda}{2}\big(\|\mb{w}\|_2^2+b^2\big).
\end{align*}
The decision variable~$\mb{w}$ represents the model weights assigned to the features and~$b$ is the bias term. It is straightforward to verify that the local loss functions~$f_i$ satisfy both Assumptions~\ref{a1:strongly_cvx} and~\ref{a2:Lipschitz}. The feature vectors,~${\mb{c}^i}^{(j)}$, are drawn from~IID Gaussian distributions with mean~$0$ and variance~$9$. We then generate binary labels from a Bernoulli distribution, with probability of~${{y^i}^{(j)}=+1}$ being~$(1+e^{\widetilde{\mb{x}}^{\top}{\mb{c}^i}^{(j)}})^{-1}$, where~$\mb{w}$ and~$b$ are drawn from the standard uniform distribution. The network topology varies according to the periodic sequence of directed graphs as shown in Fig.~\ref{fig:B_1-4} making the directed communication network~$4$-\begin{color}{black}bounded \end{color}strongly connected. 

\begin{figure}[!h]
\centering
\subfigure{\includegraphics[width=0.73\columnwidth]{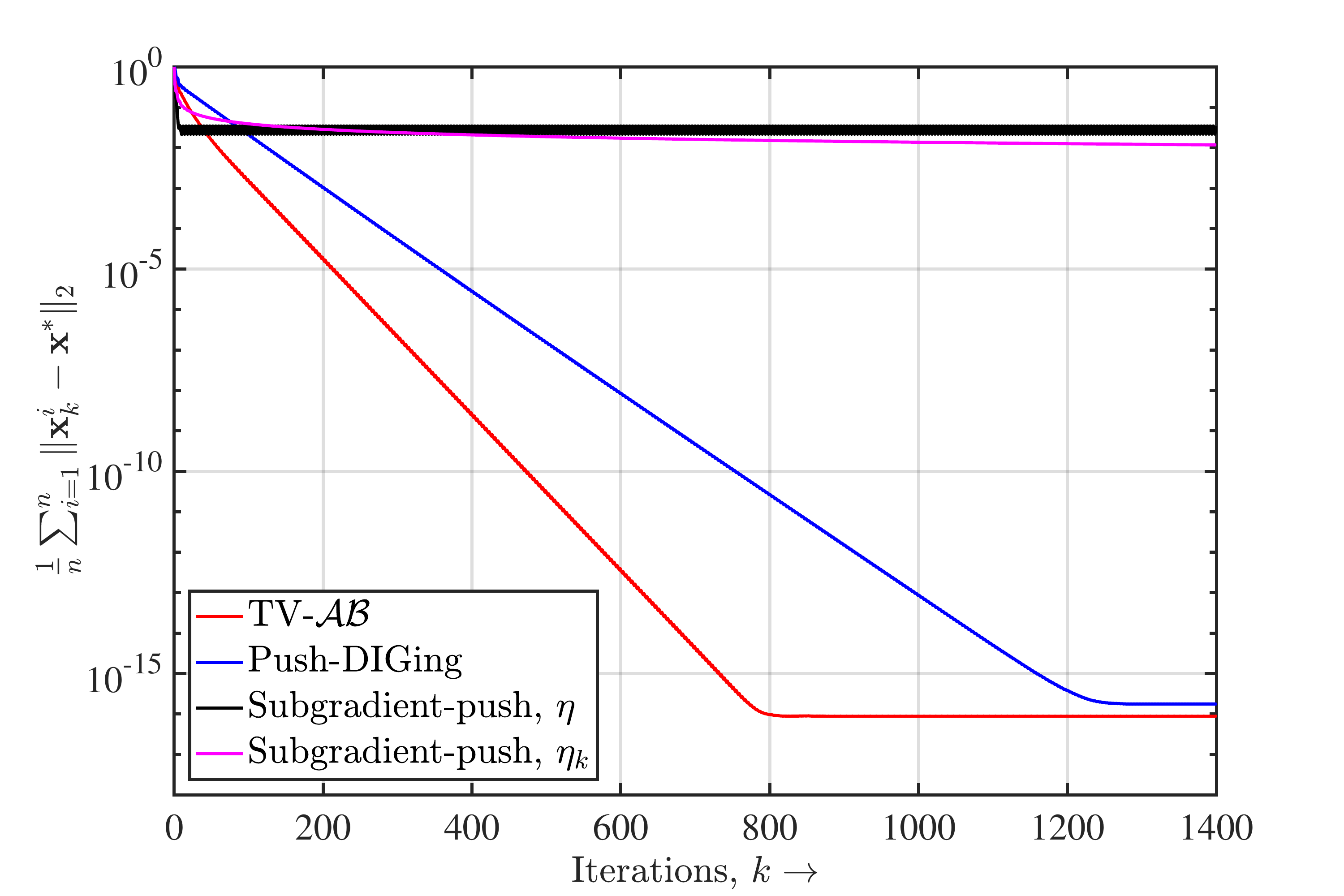}}
\subfigure{\includegraphics[width=0.73\columnwidth]{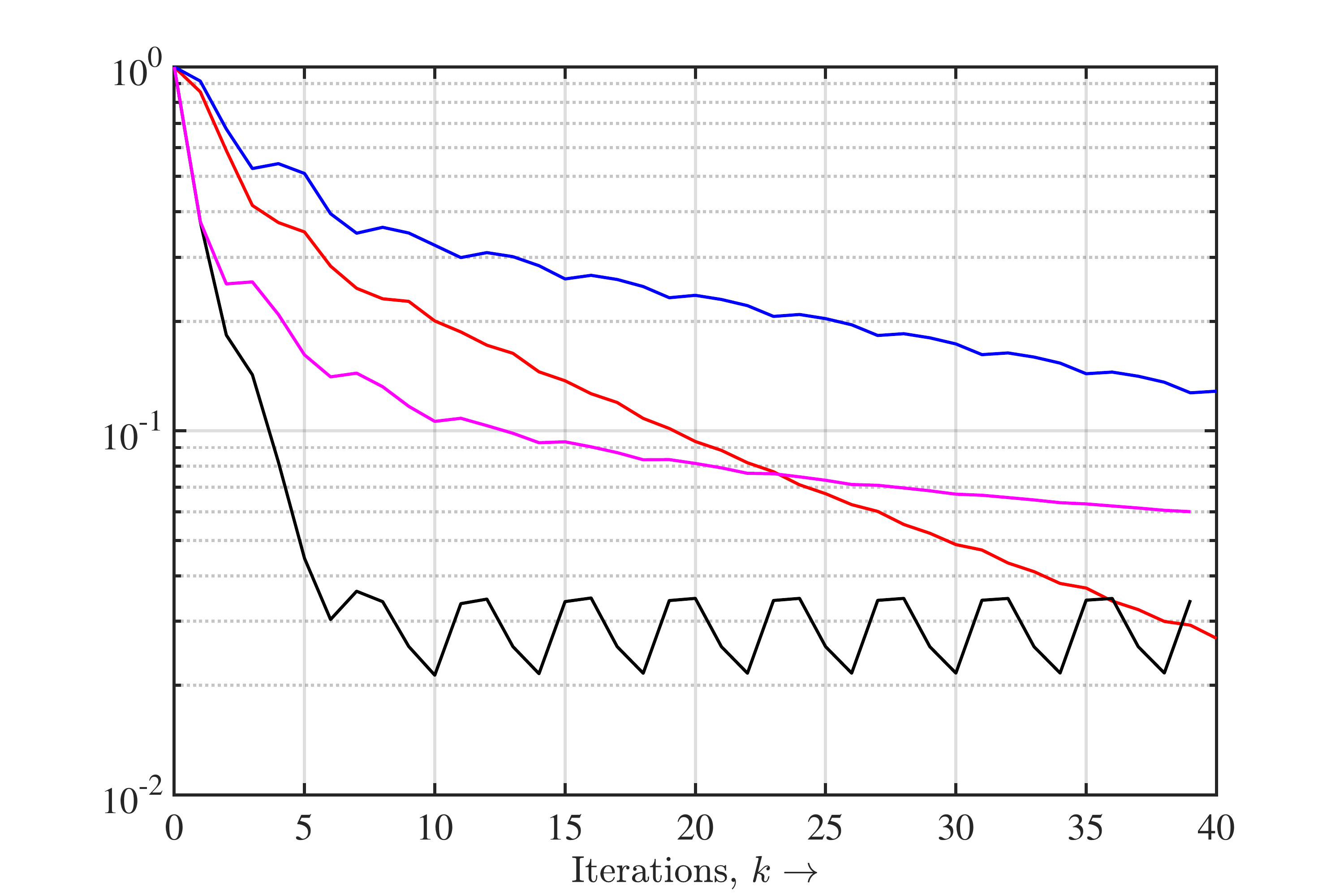}}
\caption{Distributed logistic regression: Performance comparison with the transients magnified.}
\label{fig:B_RLR_Compare2}
\end{figure}

To solve the classification problem in a distributed manner, agents initialize their states~$\mb{x}^i_0$ according to~IID zero-mean Gaussian random variables with variance~$9$. The performance of TV-$\mc{AB}$ along with Push-DIGing~\cite{diging}, and subgradient-push~\cite{opdirect_Nedic} (with constant and diminishing step-sizes) is shown in Fig.~\ref{fig:B_RLR_Compare2} with the average residual~${{1}/{n}\sum_{i=1}^n\|\mb{x}^i_k-\mb{x}^{\ast}\|_2}$ as the evaluation metric. The step-sizes for~TV-$\mc{AB}$ and Push-DIGing are hand-optimized. This numerical experiment confirms that time-varying~$\mc{AB}$ converges linearly and is observed to be faster than Push-DIGing.
\subsection{Distributed least-squares}\label{ex:D}
\begin{figure}[!ht]
\centering
\includegraphics[width=0.76\columnwidth]{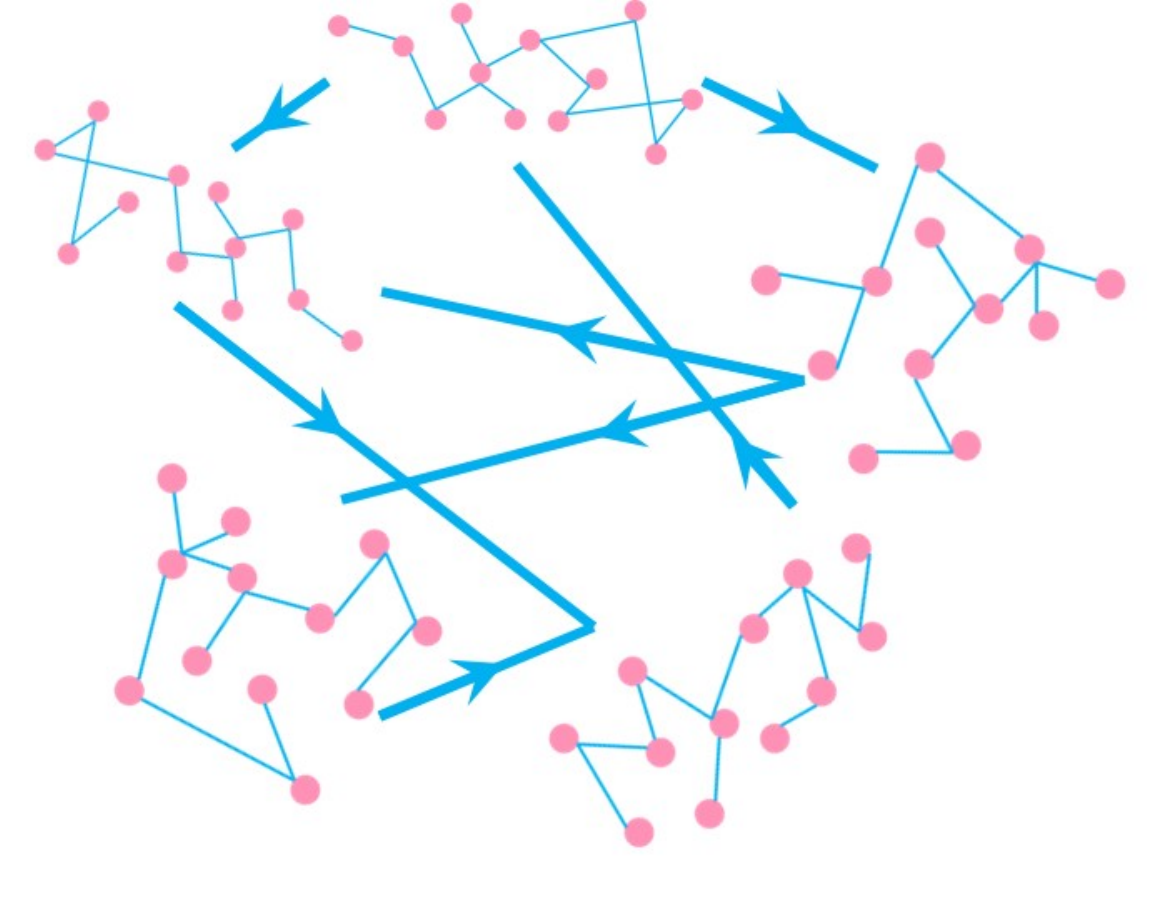}
\caption{$C=50$-bounded strongly-connected digraph: The inter-cluster graph activates every $50$th iteration.}
\label{fig:DD}
\end{figure}

In this example,~${n=60}$ agents communicate according to a time-varying sequence of~$C=50$-bounded strongly-connected digraphs. The nodes are partitioned into~$5$ equally-sized clusters, and each cluster is internally strongly-connected at every iteration. The clusters communicate with each other according to a strongly-connected cluster-level network every~$C=50$\textsuperscript{th} iteration, see~Fig.~\ref{fig:DD}. The agents aim to collaboratively find the solution~$\mb{x}^{\ast}$ of the following least-squares problem:
\[
f(\mb{x}) =\textstyle\sum_{i=1}^nf_i(\mb{x})=\frac{1}{2}\sum_{i=1}^n \|H_i\mb{x}-\mb{b}_i\|_2^2,
\]
where the vectors and matrices are of appropriate dimension. We choose each~$H_i$ such that it is rank-deficient but~$\sum_iH_i^\top H_i$ is invertible. In other words, no agent can find~$\mb{x}^*$ on its own and must cooperate. 

\begin{figure}[!h]
\centering
\subfigure{\includegraphics[width=0.73\columnwidth]{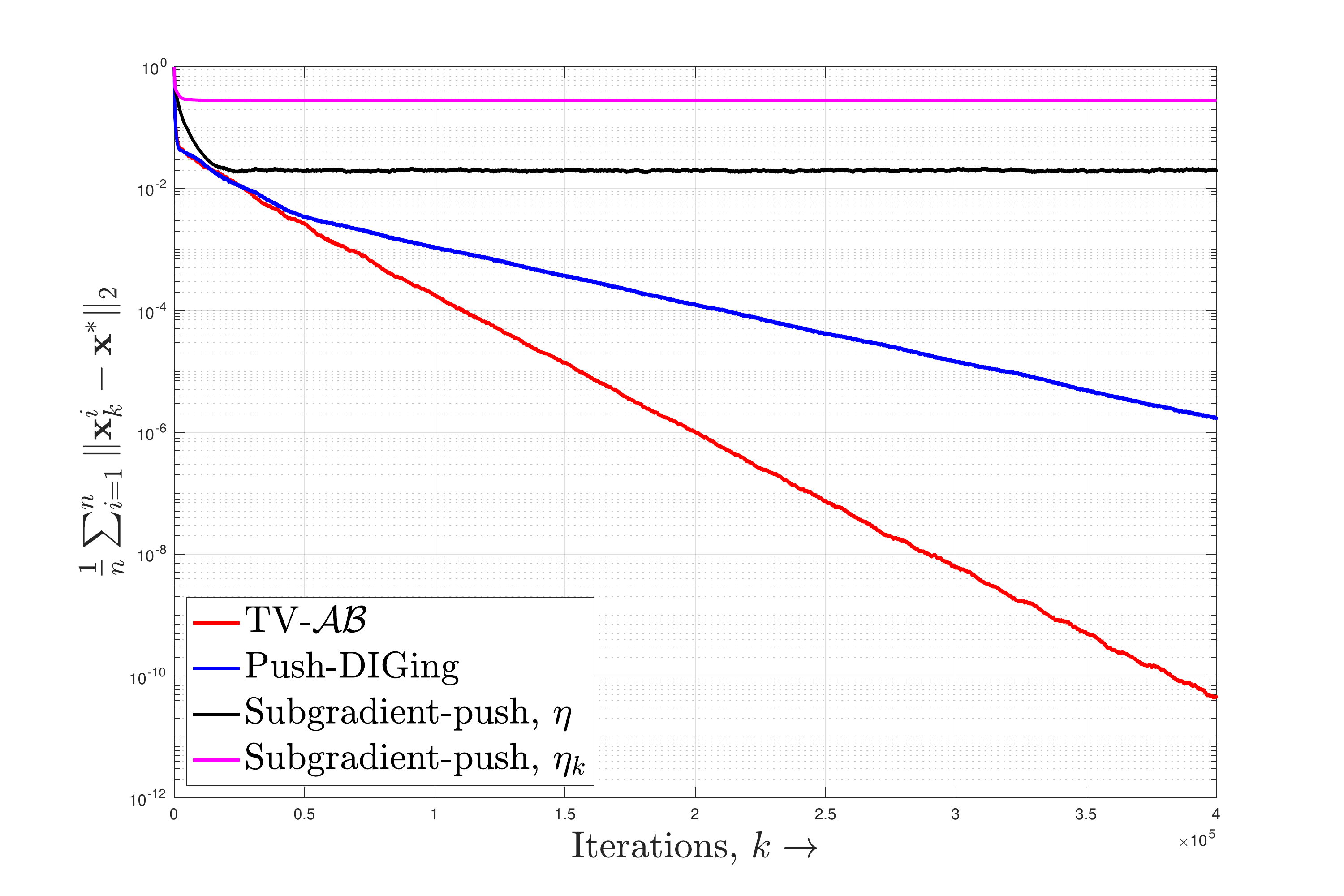}}
\subfigure{\includegraphics[width=0.73\columnwidth]{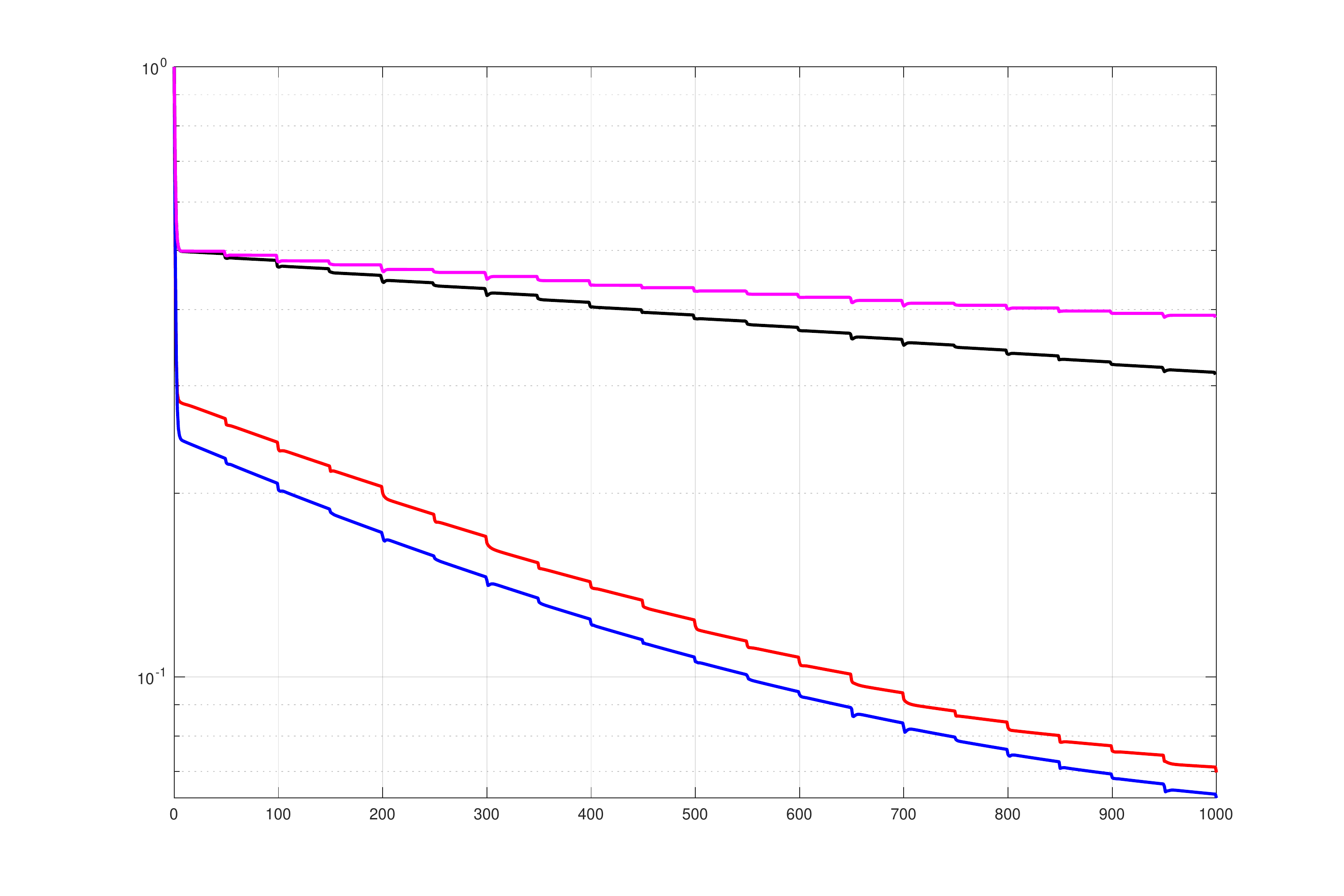}}
\caption{Distributed least squares: Performance comparison with the transients magnified.}
\label{fig:D_LS_Compare}
\end{figure}

To collaboratively solve the least-squares problem, agents initialize there states~$\mb{x}^i_0$ according to IID standard Gaussian random variables. The performance of tTV-$\mc{AB}$ is compared with Push-DIGing and subgradient-push (with both constant and diminishing step-sizes) in Fig.~\ref{fig:D_LS_Compare} with the average residual~${1}/{n}\sum_{i=1}^n\|\mb{x}^i_k-\mb{x}^{\ast}\|_2$ as the comparison metric. The step-sizes follow a similar regime as in Experiment~\ref{ex:B}. This numerical experiment, once again, confirms the linear convergence of~$\mc{AB}$ to the optimal solution provided the step-size is sufficiently small.

\subsection{TV-$\mc{AB}$ on random graphs}
\begin{figure}[!h]
\centering
\subfigure{\includegraphics[width=0.73\columnwidth]{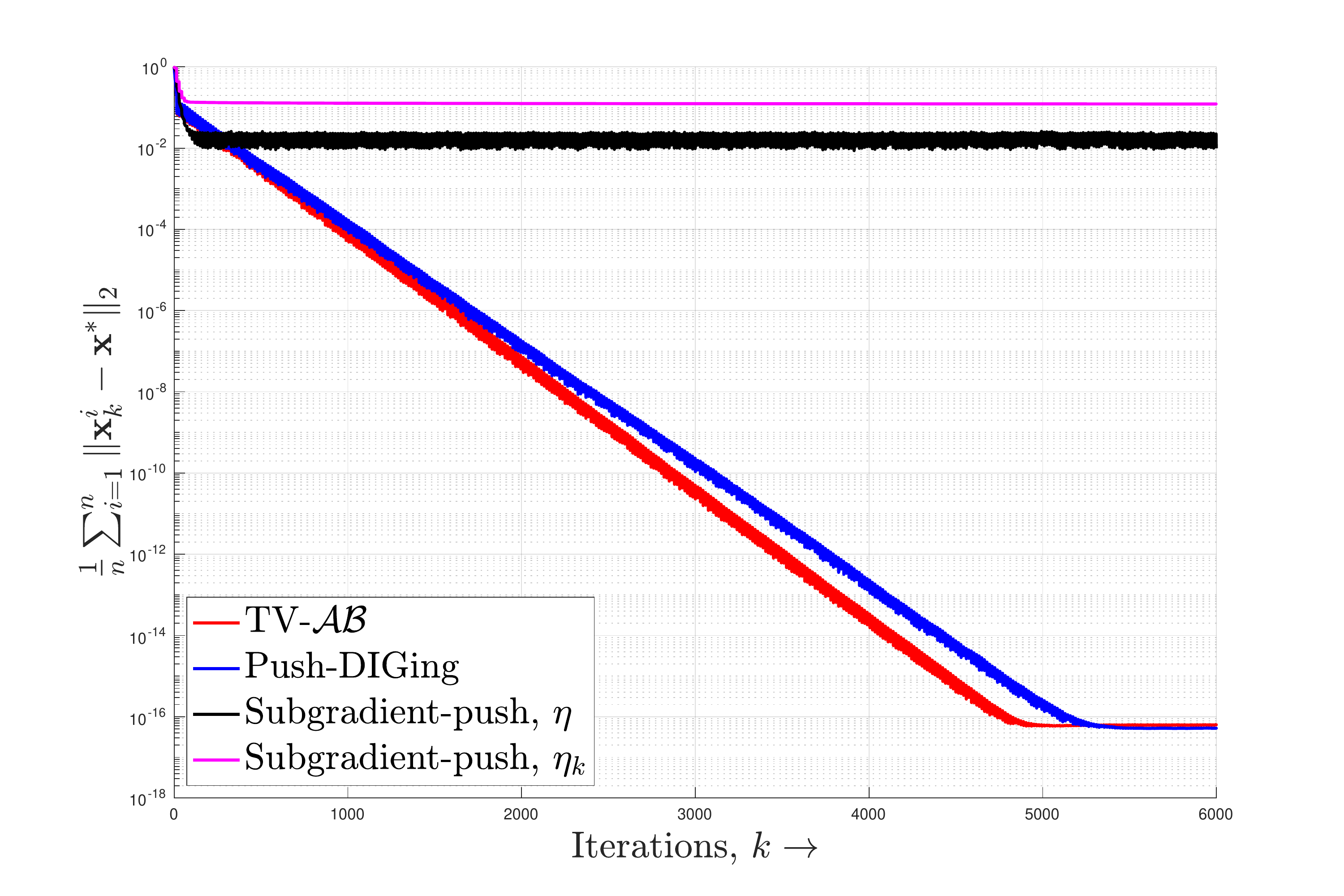}}
\subfigure{\includegraphics[width=0.73\columnwidth]{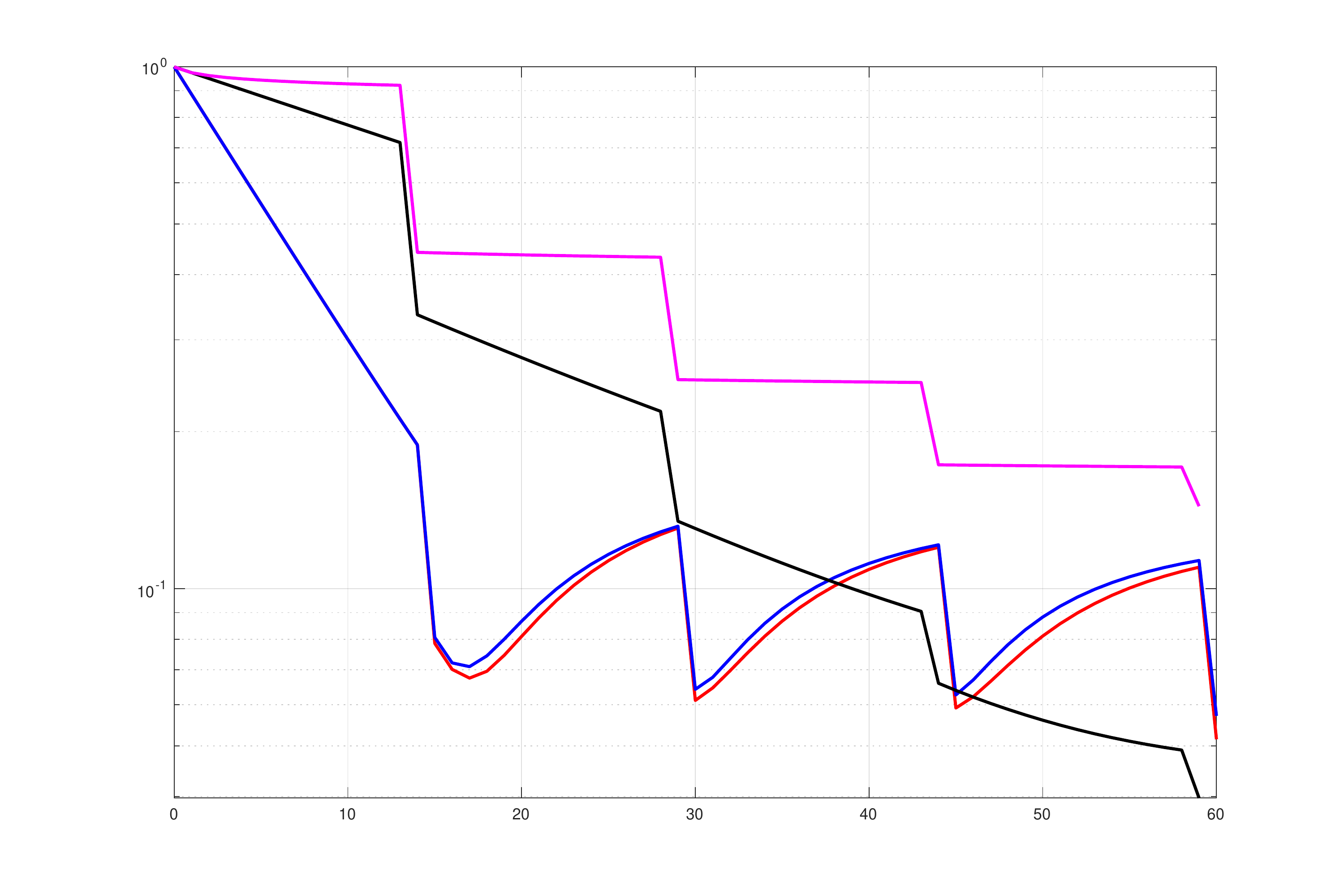}}
\caption{Distributed logistic regression over random graphs: Performance comparison  with the transients magnified.}
\label{fig:C_LogReg}
\end{figure}

We now apply~TV-$\mc{AB}$ to random networks. In this scenario,~$n=80$ agents communicate over a~$C=15$-bounded strongly-connected network to solve the logistic regression problem of example~\ref{ex:B}. The agents communicate over a strongly-connected random graph every~$15$th iteration and rely solely on local iterations for the rest of the time. The result of this experiment is presented in Fig.~\ref{fig:C_LogReg}.

\begin{figure}[!h]
\centering
\includegraphics[width=0.73\columnwidth]{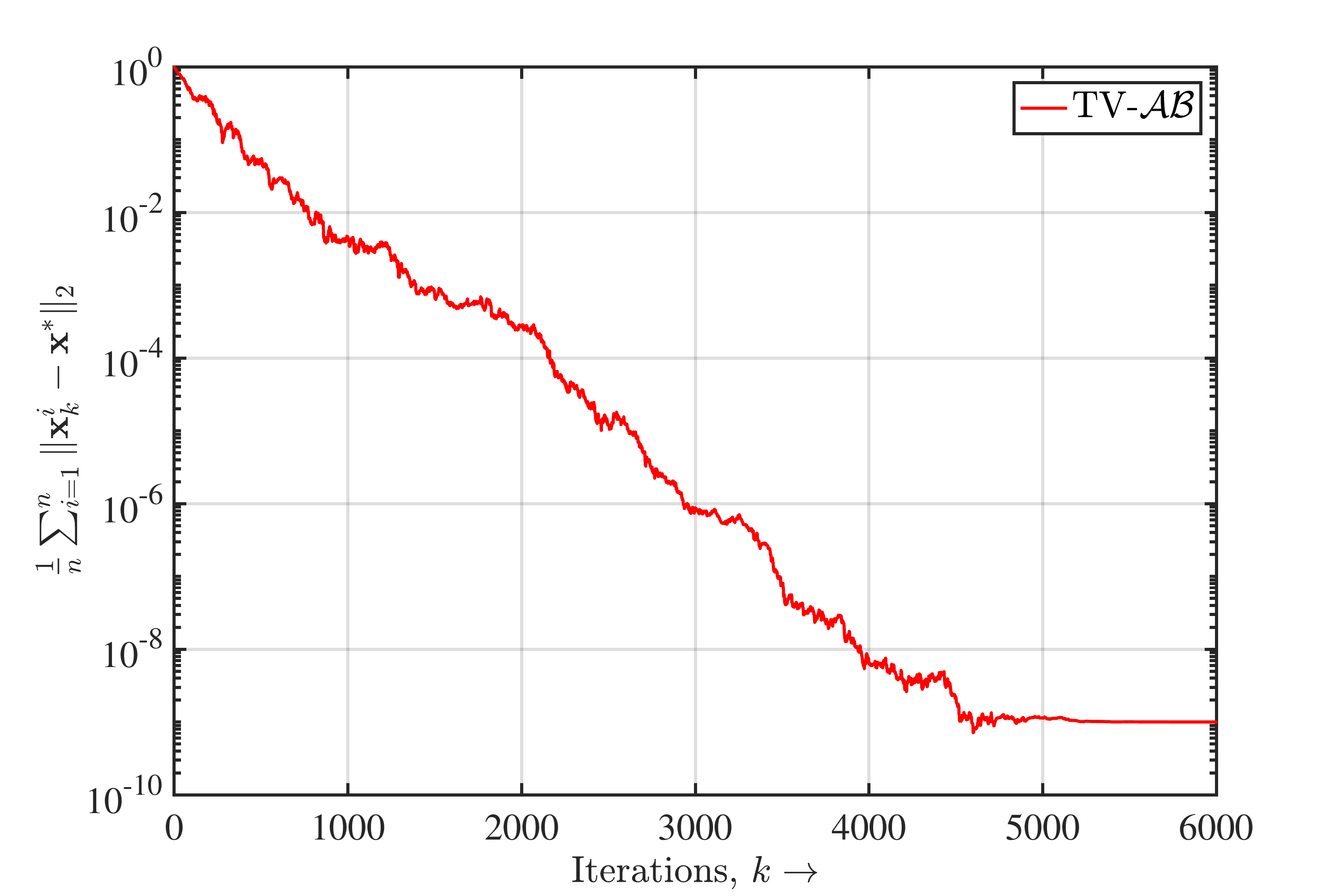}
\caption{Distributed linear regression : TV-$\mc{AB}$ under gossip.}
\label{fig:E_LinReg}
\end{figure}

In another scenario,~$10$ agents communicate according to the gossip protocol explained in Section~\ref{ss:alg}. The optimization problem is linear regression given~$10$ noisy samples of a line at each agent. The performance is compared in Fig.~\ref{fig:E_LinReg} with hand-optimized step-size.

\section{Conclusions}\label{s:conc}
In this paper, we study~TV-$\mc{AB}$ that minimizes a sum of smooth and strongly-convex functions over time-varying and possibly random directed graphs. We show that~TV-$\mc{AB}$ converges linearly to the optimal solution when underlying time-varying graphs satisfy the standard bounded connectivity assumption, i.e., a union of every~$C$ consecutive graphs is strongly-connected. We derive the convergence result based on the stability analysis of a linear system of inequalities along with a matrix perturbation argument. We further provide extensive simulations that confirm the findings in this paper. 

\appendix
\subsection{Preliminaries}\label{app:prelim}
In this section, we recap some preliminaries that are used in the analysis. We start with the definitions of absolute probability sequences~\cite{Kolmogoroff1936} and ergodicity~\cite{seneta2006non}. 
\begin{definition}[Absolute Probability Sequence] For row-stochastic matrices,~$\{R_k\}$, an absolute probability sequence is a sequence~$\{\bs{\pi}_k\}$ of stochastic vectors such that
\[\bs{\pi}_k^{\top}=\bs{\pi}_{k+1}^{\top}R_k,\quad \forall k\geq 0.\]
\end{definition}

\begin{definition}[Ergodicity] A ergodic sequence of row-stochastic matrices,~$\{R_k\}$, is such that for integers~$ p \geq 0$ and all~$i,{s = 1\,,\cdots\,,n}$
\[\lim_{c\rightarrow \infty}[U_{(c,p)}]_{i,s}\rightarrow {d}^p_s,\] 
where~$U_{(c,p)}= \mathlarger{\Pi}_{l=p}^{c}R_l$ is the backward product of~$\{R_k\}$ and~${d}^p_s$ is a constant not depending on~$i$.
\end{definition}

\noindent We next state a result on the ergodicity of the matrix sequence compliant with the aggregate digraph~${\mathcal{G}_{sC}^{{C}}= (\mathcal{V},\textstyle\cup_{l=sC}^{sC+C-1}\mathcal{E}_l)}$, see~\cite[Lemma~5.2.1]{tsitsiklisphd84} for details.
\begin{lem}
Under Assumptions~\ref{a3:joint_conn} and~\ref{a4:weight_mats}, the row-stochastic matrix sequence~${\{D_s=\textstyle\Pi_{l=sC}^{sC+C-1}A_l\}}$ compliant with the aggregate digraph~${\mathcal{G}_{sC}^{{C}}= (\mathcal{V},\cup_{l=sC}^{s+C-1}\mathcal{E}_l)}$, for all~$s\geq0$, is ergodic, i.e.,
\[
\lim_{t\rightarrow \infty}D_t\cdots D_{s+1}D_s=\mathbf{1}\bs{\mu}_s^{\top},
\]
where~$\{\bs{\mu}_s\}$ is the unique absolute probability sequence for~$\{D_s\}$ (see, e.g.,~\cite{Kolmogoroff1936},~\cite[Lemma~1]{nedic2017convergence}) and  is uniformly bounded away from zero, i.e., there exists~$\delta\in (0,\,1)$ such that~${[\bs{\mu}_s]_i\geq \delta}$,~${\forall\, i\in \mathcal{V}}$ and~${\forall\,s\geq 0}$. Furthermore, the convergence rate is geometric, i.e.,~$\forall\, t\geq s\geq 0$:
\[
\|D_t\cdots D_{s+1}D_s-\mathbf{1}\bs{\mu}_s^{\top}\|\leq \begin{color}{black}\mc{M}\end{color}q^{t-s},
\]
where the constants~$\begin{color}{black}\mc{M}\end{color}>0$ and~$q\in (0,\,1)$ depend only on~$n$ and~$\alpha$ introduced in Assumption~\ref{a4:weight_mats}.\label{lem:ergo}
\end{lem}

The next corollary extends the result above, deriving the absolute probability sequence for the sequence~$\{A_k\}$ in terms of~$\{\bs{\mu}_k\}$ in Lemma~\ref{lem:ergo}.
\begin{cor}\label{cor:phi_cor}
Under the assumptions of Lemma~\ref{lem:ergo}, the sequence~$\{\bs{\phi}_k\}$ is an absolute probability sequence for the matrix sequence~$\{A_k\}$, where
{\small \begin{align}
\begin{aligned}
\bs{\phi}_k^{\top}=&\bs{\mu}_k^{\top},\quad k=sC,\\
\bs{\phi}_k^{\top}=&\bs{\mu}_{(s+1)C}^{\top} A_{(s+1)C-1}\cdots A_k,\quad k\in\big(sC,(s+1)C\big),
\end{aligned}\label{eq:phi_t}
\end{align}}
for some~$s\geq 0$. 
\end{cor}
\begin{IEEEproof}
Since the products~$ A_{(s+1)C-1}\cdots A_k$ are row-stochastic and~$\{\bs{\mu}_s\}$ is an absolute probability sequence from Lemma~\ref{lem:ergo}, each~$\bs{\mu}_s$ is a stochastic vector by definition and so is the vector~$\bs{\phi}_k$ in Eq.~\eqref{eq:phi_t}. In fact,
\begin{align*}
\bs{\phi}_k^{\top}=&\bs{\mu}_{(s+1)C}^{\top} A_{(s+1)C-1}\cdots A_k,\\
=& \bs{\mu}_{(s+1)C}^{\top} A_{(s+1)C-1}\cdots A_{k+1} ~ A_k=\bs{\phi}_{k+1}^{\top}A_k,
\end{align*}
i.e., the sequence~$\{\bs{\phi}_k\}$ is an absolute probability sequence for~$\{A_k\}$. 
\end{IEEEproof}

The next lemma from~\cite[Lemma~5.3]{diging} establishes multi-step contraction of a backward product of a series of row-stochastic matrices~$\{A_k\}$. This result is fundamental to the convergence analysis of TV-$\mc{AB}$.
\begin{lem}[$\overline{C}_{\mc{A}}$-step contraction for~$\{\mc{A}_k\}$]Let Assumptions~\ref{a3:joint_conn} and~\ref{a4:weight_mats} hold. Recall that~$\mc{A}_k=A_k\otimes I_p$ and define an integer $\overline{C}_{\mc{A}}\geq C$ such that
\begin{align}
\gamma_\mc{A}\triangleq Q_{\mc{A}}(1-\alpha^{nC})^{\frac{\overline{C}_\mc{A}-1}{nC}}<1,\quad Q_{\mc{A}}=2n\frac{1+\alpha^{-nC}}{1-\alpha^{nC}}.
\end{align}
Then for any~$k\geq\overline{C}_{\mc{A}}-1$ and any vector~$\mb{b}\in\mbb{R}^{np}$, if 
\[\mb{a}={\mc{A}_{(\overline{C}_{\mc{A}},k-(\overline{C}_{\mc{A}}-1))}}\mb{b},\] 
where~$\mc{A}_{(\overline{C}_{\mc{A}},k-(\overline{C}_{\mc{A}}-1))}\triangleq{\mc{A}_k\cdots \mc{A}_{k-(\overline{C}_{\mc{A}}-1)}}$, we have 
\begin{align*}
&\big\|\big((I_{n}-\mb{1}_{n}\bs{\phi}_{k+1}^{\top})\otimes I_p\big)\mb{a}\big\|_2\\
&\leq \gamma_{\mc{A}} \big\|\big((I_{n}-\mb{1}_{n}\bs{\phi}_{k-(\overline{C}_{\mc{A}}-1)}^{\top})\otimes I_p\big)\mb{b}\big\|_2,
\end{align*}
where~$\{\bs{\phi}_k\}$ is the absolute probability sequence of~$\{A_k\}$, defined in~Eq.~\eqref{eq:phi_t}.\label{lem:contract_A}
\end{lem}
The next corollary establishes the multi-step contraction for the sequence~${\{R_k\}={V}^{-1}_{k+1}{B}_k{V}_k}$, where~${V_k=\textrm{diag}[\mb{v}_k]}$ and~$\{\mb{v}_k\}$ evolves as Eq.~\eqref{eq:left_eig}.

\begin{cor}[{$\overline{C}_{\mc{B}}$-step contraction for~$\{\mc{R}_k\}$}]
Let Assumptions~\ref{a3:joint_conn} and~\ref{a4:weight_mats} hold. Recall that~${\mc{R}_k=R_k\otimes I_p}$ and define an integer~${\overline{C}_{\mc{B}}\geq C}$ such that
\begin{align}
\gamma_{\mc{B}}\triangleq Q_{\mc{B}}(1-\tau^{nC})^{\frac{\overline{C}_{\mc{B}}-1}{nC}}<1,\quad& Q_{\mc{B}}=2n\frac{1+\tau^{-nC}}{1-\tau^{nC}},
\end{align}
where~$\tau=\frac{\beta}{n^{nC+1}}$; then for any~$k\geq\overline{C}_{\mc{B}}-1$ and any vector~${\mb{b}\in\mbb{R}^{np}}$, if 
\[\mb{a}={\mc{R}_{(\overline{C}_{\mc{B}},k-(\overline{C}_{\mc{B}}-1))}}\mb{b},\] 
where~${\mc{R}_{(\overline{C}_{\mc{B}},k-(\overline{C}_{\mc{B}}-1))}}={\mc{R}_k\cdots \mc{R}_{k-(\overline{C}_{\mc{B}}-1)}}$, we have 
\begin{align*}
&\big\|\big((I_{n}-\mb{1}_{n}\mb{v}_{k+1}^{\top})\otimes I_p\big)\mb{a}\big\|_2\\
\leq &\gamma_{\mc{B}} \big\|\big((I_{np}-\mb{1}_{n}\mb{v}_{k-(\overline{C}_{\mc{B}}-1)}^{\top})\otimes I_p\big)\mb{b}\big\|_2.
\end{align*}
\label{cor:contract_B}\end{cor}
\begin{IEEEproof} 
It can be verified that
{\small \begin{align*}
&((I_{n}-\mb{1}_n\mb{v}_{k+1}^{\top})\otimes I_p)\mc{R}_{(\overline{C}_{\mc{B}},k-(\overline{C}_{\mc{B}}-1))}\\
=&(({R}_{(\overline{C}_{\mc{B}},k-(\overline{C}_{\mc{B}}-1))}-\mb{1}_n\mb{v}_{k-(\overline{C}_{\mc{B}}-1)}^{\top})(I_n-\mb{1}_n\mb{v}_{k-(\overline{C}_{\mc{B}}-1)}^{\top}))\otimes I_p.
\end{align*}}Therefore,
{\small \begin{align*}
&\big\|\big((I_{n}-\mb{1}_n\mb{v}_{k+1}^{\top})\otimes I_p\big)\mb{a}\big\|_2\\
\leq&\big\|{R}_{(\overline{C}_{\mc{B}},k-(\overline{C}_{\mc{B}}-1))}-\mb{1}_n\mb{v}_{k-(\overline{C}_{\mc{B}}-1)}^{\top}\big\|_2\cdot\\
&\big\|((I_{n}-\mb{1}_n\mb{v}_{k-(\overline{C}_{\mc{B}}-1)}^{\top})\otimes I_p)\mb{b}\big\|_2,
\end{align*}}where the inequality follows from the compatibility of vector~$2$-norm with matrix spectral norm. We now find an upper bound for the first term on the right hand side of the above equation.
First note that~${[R_{k}]_{i,j} = [B_k]_{ij}{[\mb{v}_k]_j}/{[\mb{v}_{k+1}]_i}}$. From~{\cite[Corollary~2(b)]{opdirect_Nedic}}, we have that~${[\mb{v}_k]_j\geq {1}/{n^{nC}}}$,~${\forall\,k\geq 0}$. Since~${1}/{[\mb{v}_k]_j}\geq{1}/{n}$ and for any~$(i,j)\in\mathcal{E}_k$,~${[B_{k}]_{i,j}\geq \beta}$ by Assumption~\ref{a4:weight_mats}, we have that in~${[R_{k}]_{i,j}\geq \tau\triangleq {\beta}/{n^{nC+1}}}$, for any~$(i,j)\in\mathcal{E}_k$. Therefore, noting that for an~${n\times n}$ matrix,~$X$,~${\|X\|_2\leq n\|X\|_{\max}}$, we have
\begin{align*}
&\big\|{R}_{(\overline{C}_{\mc{B}},k-(\overline{C}_{\mc{B}}-1))}-\mb{1}_n\mb{v}^{\top}_{k-(\overline{C}_{\mc{B}}-1)}\big\|_2 \\
\leq& n\big\| R_{(\overline{C}_B,k-(\overline{C}_B-1))}-\mb{1}_n\mb{v}^{\top}_{k-(\overline{C}_B-1)}\big\|_{\max}\leq\gamma_{\mc{B}},
\end{align*}
where~$\gamma_{\mc{B}}\triangleq 2n \frac{1+\tau^{-nC}}{1-\tau^{nC}}(1-\tau^{nC})^{\frac{\overline{C}_{\mc{B}}-1}{nC}}$ and
the last inequality is from~\cite[Lemma~4(c)]{uc_Nedic}.
\end{IEEEproof}

Finally, the next lemma is a standard result in optimization theory and states that the optimality gap in the domain space shrinks by at least a fixed ratio for a  gradient descent step.
\begin{lem}[\begin{color}{black}{\cite[Lemma~3.11]{opt_literature0}}\end{color}] Let~$g:\mathbb{R}^p\mapsto\mathbb{R}$ be~$\mu$-strongly-convex and have~$\ell$-Lipschitz gradient. Define~${\mb{x}^+ = \mb{x}-\zeta {\nabla}g(\mb{x})}$, where~$0 <\zeta< {2}/{\ell}$. Then
\[\|\mb{x}^+-\mb{x}^{\ast}\|_2\leq  \chi\|\mb{x}-\mb{x}^{\ast}\|_2\]
where~$\chi= \max\{|1-\zeta\mu|,|1-\zeta \ell|\}$.\label{lem:har}
\end{lem}

In the following, we provide the proofs of the lemmas described earlier in the paper. 
\subsection{Proof of Lemma~\ref{lem:y}.}\label{app:lem_4}
\begin{IEEEproof}
Recall that~$\mathbf{s}_k=\widetilde{\mathbf{s}}^{\textrm{w}}_k+(\mb{1}_n\mb{v}_k^{\top}\otimes I_p){\mathbf{s}}_k$ and it can be verified that
\begin{align}
(\mb{1}_n\mb{v}_k^{\top}\otimes I_p){\mathbf{s}}_k=(\mb{1}_n\mb{1}_n^{\top}\otimes I_p){\nabla \mb{f}}(\mathbf{x}_k).
\label{eq:11}
\end{align}
Exploiting the optimality condition~${{\nabla \mb{f}}(\mb{1}_n\otimes \mathbf{x}^{\ast})=\mb{0}_{np}}$, we can express~$\mathbf{s}_k$ as
\begin{align*}
\mathbf{s}_k=&\widetilde{\mathbf{s}}^{\textrm{w}}_k+(\mb{1}_n\mb{1}_n^{\top}\otimes I_p)({\nabla \mb{f}}(\mathbf{x}_k)-{{\nabla \mb{f}}(\mb{1}_n\otimes \mathbf{x}^{\ast})}).
\end{align*}
Therefore, from the triangle inequality we have
{\small\begin{align*}
\|\mathbf{s}_k\|_2 & \leq \| \widetilde{\mathbf{s}}^{\textrm{w}}_k\|_2+\sqrt{n}\textstyle\sum_{i=1}^n\| {\nabla \mb{f}}_i(\mb{x}^i_k)-{\nabla \mb{f}}_i(\mb{x}^{\ast})\|_2,\\
&\leq \| \widetilde{\mathbf{s}}^{\textrm{w}}_k\|_2+\sqrt{n}L\textstyle\sum_{i=1}^n\| \mb{x}^i_k-\mb{x}^{\ast}\|_2,\\
&\leq \| \widetilde{\mathbf{s}}^{\textrm{w}}_k\|_2+nL\big\|\mb{x}_k+(\mb{1}_n\phi^{\top}_k\otimes I_p)(\mb{x}_k-\mb{x}_k)-\mb{1}_n\otimes \mb{x}^{\ast}\big\|_2,\\
&\leq \| \widetilde{\mathbf{s}}^{\textrm{w}}_k\|_2+{n}L\| \widetilde{\mathbf{x}}^{\textrm{w}}_k\|_2+{n}L\| {\mathbf{r}}_k\|_2.
\end{align*}}where the second inequality uses Lipschitz continuity and the third inequality is a consequence of Cauchy-Schwarz. Noting that~$\mathbf{y}_k=(V_k\otimes I_p)\mb{s}_k$, we have $\|\mathbf{y}_k\|_2{\leq}\|(V_k\otimes I_p)\|_2\|\mb{s}_k\|_2$ from the compatibility of vector~$2$-norm with matrix spectral norm. Next, since~${\|(V_k\otimes I_p)\|_2=\|V_k\|_2=\max_i [\mb{v}_k]_i}<1$, we have
\begin{align*}
\|\mathbf{y}_k\|_2\leq&\|\mb{s}_k\|_2\leq{n}L\| \widetilde{\mathbf{x}}^{\textrm{w}}_k\|_2+{n}L\| {\mathbf{r}}_k\|_2+\| \widetilde{\mathbf{s}}^{\textrm{w}}_k\|_2,
\end{align*}
and the lemma follows.
\end{IEEEproof}

\subsection{Proof of Lemma~\ref{lem:x_tilde}.}\label{app:lem_5}
\begin{IEEEproof}
From Eq.~\eqref{eq:state}, we have
\begin{align*}
\mathbf{x}_{k+1}=&\mc{A}_{(\overline{C},k-(\overline{C}-1))}\mathbf{x}_{k-(\overline{C}-1)}-\eta\sum_{l=0}^{\overline{C}-1} \mc{A}_{l,k-(l-1)}{\mb{y}_{k-l}},
\end{align*}
which leads to 
{\small\begin{align*}
\|\widetilde{\mathbf{x}}^{\textrm{w}}_{k+1}\|_2=&\big\|\big((I_n-\mb{1}_n\bs{\phi}^{\top}_{k+1})\otimes I_p\big)\mathbf{x}_{k+1}\big\|_2\\
{\leq}&\big\|\big((I_n-\mb{1}_n\bs{\phi}^{\top}_{k+1})\otimes I_p\big)\mc{A}_{(\overline{C},k-(\overline{C}-1))}\mathbf{x}_{k-(\overline{C}-1)}\big\|_2\\
&+\eta\textstyle\sum_{l=0}^{\overline{C}-1} \big\|\big((I-\mb{1}_n\bs{\phi}^{\top}_{k+1})\otimes I_p\big)\mc{A}_{l,k-(l-1)}\mathbf{y}_{k-l}\big\|_2.
\end{align*}}Consequently,~$\forall\, k\geq \overline{C}-1$,
\begin{align*}
\|\widetilde{\mathbf{x}}^{\textrm{w}}_{k+1}\|_2
\leq & \gamma_{\mc{A}}\|\widetilde{\mathbf{x}}^{\textrm{w}}_{k-(\overline{C}-1)}\|_2+\eta Q_{\mc{A}}\textstyle\sum_{l=0}^{\overline{C}-1}\|\mathbf{y}_{k-l}\|_2,
\end{align*}
where the second term follows from~{\cite[Lemma~4(c)]{uc_Nedic}}. From Lemma~\ref{lem:y}, we have
{\small\begin{align*}
 \|\widetilde{\mathbf{x}}^{\textrm{w}}_{k+1}\|_2
\leq& ({\gamma_{\mc{A}}}+\eta Q_{\mc{A}}{n}L)\|\widetilde{\mathbf{x}}^{\textrm{w}}_{k-(\bar{C}-1)}\|_2,\\
&+\eta Q_{\mc{A}}{n}L\textstyle\sum_{l=0}^{\overline{C}-2}\|\widetilde{\mathbf{x}}^{\textrm{w}}_{k-l}\|_2,\\
&+\eta Q_{\mc{A}}{n}L\|{\mathbf{r}}_{k-(\bar{C}-1)}\|_2+\eta Q_{\mc{A}}{n}L\textstyle\sum_{l=0}^{\overline{C}-2}\|{\mathbf{r}}_{k-l}\|_2,\\
&+\eta Q_{\mc{A}}\|\widetilde{\mathbf{s}}^{\textrm{w}}_{k-(\bar{C}-1)}\|_2+\eta Q_{\mc{A}}\textstyle\sum_{l=0}^{\overline{C}-2}\|\widetilde{\mathbf{s}}^{\textrm{w}}_{k-l}\|_2,
\end{align*}}where~$\gamma_{\mc{A}}$ and~$Q_{\mc{A}}$ are the constants defined in Lemma~\ref{lem:contract_A}.
\end{IEEEproof}

\subsection{Proof of Lemma~\ref{lem:r}.}\label{app:lem_6}
\begin{IEEEproof}
Note that~$\mathbf{y}_k-(\mb{v}_k\mb{1}_n^{\top}\otimes I_p)\mathbf{y}_k={(V_k\otimes I_p)\widetilde{\mathbf{s}}^{\textrm{w}}_k}$, we~have
{\small\begin{align}
\|{\mathbf{r}}_{k+1}\|_2&=\big\|(\mb{1}_n\bs{\phi}_{k+1}^{\top}\otimes I_p)\mathbf{x}_{k+1}-\mathbf{1}_n\otimes {\mb{x}^{\ast}}\big\|_2\nonumber\\
&=\big\|(\mb{1}_n\bs{\phi}_{k+1}^{\top}\otimes I_p)\big(\mc{A}_k\mathbf{x}_k-\eta\mathbf{y}_k\nonumber\\
&+(\mb{v}_k\mb{1}_n^{\top}\otimes I_p)\mathbf{y}_k(-\eta+\eta)\big)-\mathbf{1}_n\otimes {\mb{x}^{\ast}}\big\|_2\nonumber\\
&{\leq} \eta \big\| (\mb{1}_n\bs{\phi}_{k+1}^{\top}\otimes I_p)\big(\mathbf{y}_k-(\mb{v}_k\mb{1}_n^{\top}\otimes I_p)\mathbf{y}_k\big)\big\|_2\nonumber\\
&+\big\|(\mb{1}_n\bs{\phi}_k^{\top}\otimes I_p)\mathbf{x}_k-\mathbf{1}_n\otimes {\mb{x}^{\ast}}-\eta
\theta_k(\mb{1}_n\mb{1}_n^{\top}\otimes I_p)\mathbf{y}_k\big\|_2\nonumber\\
&\leq \eta \sqrt n\|\widetilde{\mathbf{s}}^{\textrm{w}}_k\|_2+\big\|(\mb{1}_n\bs{\phi}_k^{\top}\otimes I_p)\mathbf{x}_k-\mathbf{1}_n\otimes {\mb{x}^{\ast}}\nonumber\\
&\,\,\quad\qquad\qquad\qquad-\eta \theta_k(\mb{1}_n\mb{1}_n^{\top}\otimes I_p)\nabla\mathbf{f}(\mb{x}_k)\big\|_2,\label{eq:go_back}
\end{align}}where~$\theta_k=\bs{\phi}_{k+1}^{\top}\mb{v}_k$. From~\cite[Lemma~4(c)]{uc_Nedic}, we note that~${\theta_k\geq{1}/{n^{nC}}}$, while~$\theta_k\leq1,\forall k,$ from Cauchy-Schwarz. The substitution of~$(\mb{1}_n\mb{1}_n^{\top}\otimes I_p)\mathbf{y}_k$ by~$(\mb{1}_n\mb{1}_n^{\top}\otimes I_p)\nabla \mathbf{f}(\mb{x}_k)$ follows a similar reasoning as in Eq.~\eqref{eq:11}. Furthermore, the second term on the right hand side of the above inequality can be expressed as follows:
\begin{align*}
&\big\|{\mb{1}_n\otimes \overline{\mb{x}}^{\textrm{w}}_k}-n\eta\theta_k\nabla\mb{f}(\mb{1}_n\otimes \overline{\mb{x}}^{\textrm{w}}_k)-\mathbf{1}_n\otimes{\mb{x}^{\ast\top}}\\
&+n\eta\theta_k\nabla\mb{f}(\mb{1}_n\otimes \overline{\mb{x}}^{\textrm{w}}_k)
-\eta \theta_k(\mb{1}_n\mb{1}_n^{\top}\otimes I_p){\nabla \mb{f}}(\mathbf{x}_k)\big\|_2\\
{\leq} & \big\| \mathbf{1}_n\otimes \overline{\mb{x}}^{\textrm{w}}_k-n\eta\theta_k{\nabla \mb{f}}(\mb{1}_n\otimes \overline{\mb{x}}^{\textrm{w}}_k)-\mathbf{1}_n\otimes{\mb{x}^{\ast\top}}\big\|_2\\
&+\eta\theta_k\big\|\mb{1}_n\otimes\big(\textstyle\sum_{i=1}^n({\nabla f}_i(\overline{\mb{x}}^{\textrm{w}}_k)-{\nabla f}_i({\mb{x}}^i_k))\big)\big\|_2\\
{\leq} & \sqrt n\| \overline{{\mb{x}}}^{\textrm{w}}_k-n\eta\theta_k{\nabla f}(\overline{\mb{x}}^{\textrm{w}}_k)-{\mb{x}^{\ast}}\|_2+\eta{n}L\|\overline{\mathbf{x}}^{\textrm{w}}_k-\mb{x}_k\|_2,
\end{align*}
where the second term is due to Assumption~\ref{a2:Lipschitz}. From Lemma~\ref{lem:har}, if~$0<\eta<\frac{2}{nL}<\frac{2}{nL\theta_k}$, we can bound the first term on the right hand side and obtain
{\small\begin{align*}
\sqrt{n}\| \overline{{\mb{x}}}^{\textrm{w}}_k-n\eta\theta_k{\nabla f}(\overline{\mb{x}}^{\textrm{w}}_k)-{\mb{x}^{\ast}}\|_2 \leq & \sqrt n\chi_k \| \overline{{\mb{x}}}^{\textrm{w}}_k-{\mb{x}^{\ast}}\|_2=\chi_k \| \mb{r}_k\|_2,
\end{align*}}where 
\begin{align*}
\chi_k=&\max\{|1-n\eta\theta_k{\mu}|,|1-n\eta\theta_k L|\}\\
=&1-n\eta\theta(k){\mu}\leq 1-\eta\frac{{\mu}}{n^{nC-1}}
\end{align*} 
for~${\eta<\frac{1}{nL}<\frac{1}{nL\theta_k}}$.
Going back to Eq.~\eqref{eq:go_back},
\begin{align*}
\|{\mathbf{r}}_{k+1}\|_2\leq & \eta \sqrt{n}\|\widetilde{\mathbf{s}}^{\textrm{w}}_k\|_2+\chi_k \|\mb{r}_k\|_2+\eta n L {\|\widetilde{\mathbf{x}}^{\textrm{w}}_k\|_2}\\
\leq &\eta nL\|\widetilde{\mathbf{x}}^{\textrm{w}}_k\|_2+\big(1-\eta\frac{{\mu}}{n^{nC-1}}\big)\|\mb{r}_k\|_2+\eta\sqrt{n} \|\widetilde{\mathbf{s}}^{\textrm{w}}_k\|_2,
\end{align*}
and the lemma follows.
\end{IEEEproof}

\subsection{Proof of Lemma~\ref{lem:s_tilde}.}\label{app:lem_7}
\begin{IEEEproof}
From Eq.~\eqref{eq:s_tracking}, we have
\begin{align*}
\mathbf{s}_{k+1}=&\mc{R}_{(\overline{C},k-(\overline{C}-1))}\mathbf{s}_{k-(\overline{C}-1)}\nonumber\\
&+\sum_{l=0}^{\overline{C}-1} \mc{R}_{(l,k-(l-1))}(V^{-1}_{k-(l-1)}\otimes I_p)\mathbf{z}_{k-(l-1)}.
\end{align*}
Applying triangle inequality, we get
{\small \begin{align*}
&\|\widetilde{\mathbf{s}}^{\textrm{w}}_{k+1}\|_2=\|((I_n-\mb{1}_n\mb{v}_{k+1}^{\top})\otimes I_p)\mathbf{s}_{k+1}\|_2\\
{\leq}&\|((I_n-\mb{1}_n\mb{v}_{k+1}^{\top})\otimes I_p)\mc{R}_{(\overline{C},k-(\overline{C}-1))}\mathbf{s}_{k-(\overline{C}-1)}\|_2\\
+&\sum_{l=-1}^{\overline{C}-2} \|((I_n-\mb{1}_n\mb{v}_{k+1}^{\top})\otimes I_p)\mc{R}_{(l+1,k-l)}(V^{-1}_{k-l}\otimes I_p)\mathbf{z}_{k-l}\|_2,
\end{align*}}where~$\mc{R}_{(l,k-l)}=I_{np}$ for~$l\leq0$. Therefore,~$\forall\, k\geq \overline{C}-1$,
\begin{align*}
\|\widetilde{\mathbf{s}}^{\textrm{w}}_{k+1}\|_2\leq &{\gamma_{\mc{B}}}\|((I_n-\mb{1}_n\mb{v}_{k-(\overline{C}-1)}^{\top})\otimes I_p)\mathbf{s}_{k-(\overline{C}-1)}\|_2\\
&+{Q_{\mc{B}}}\textstyle\sum_{l=-1}^{\overline{C}-2} \|(V^{-1}_{k-l}\otimes I_p)\mathbf{z}_{k-l}\|_2,
\end{align*}
where the second term follows~\cite[Lemma~4(c)]{uc_Nedic}. With the help of~{\cite[Corollary~2(b)]{opdirect_Nedic}} and Corollary~\ref{cor:contract_B}, we have
\begin{align*}
\|\widetilde{\mathbf{s}}^{\textrm{w}}_{k+1}\|_2\leq &\gamma_{\mc{B}}\|\widetilde{\mathbf{s}}^{\textrm{w}}_{k-(\overline{C}-1)}\|_2+n^{nC}Q_{\mc{B}}\textstyle\sum_{l=0}^{\overline{C}-1} \|\mathbf{z}_{k-(l-1)}\|_2.
\end{align*}
From Assumption~\ref{a2:Lipschitz}, the summation in the second term can be bounded as follows:
\begin{align*}
&\textstyle\sum_{l=0}^{\overline{C}-1} \|\mathbf{z}_{k-(l-1)}\|_2 {\leq} L\sum_{l=0}^{\overline{C}-1}\|\mathbf{x}_{k-(l-1)}-\mathbf{x}_{k-l}\|_2.
\end{align*}
Furthermore,
\begin{align*}
& \| \mb{x}_{k-(l-1)}-\mb{x}_{k-l}\|_2\\
= & \big\| (\mc{A}_{k-l}-I_{np}) \big(\mb{x}_{k-l}-(\mb{1}_n\bs{\phi}^{\top}_{k-l}\otimes I_p )\mb{x}_{k-l}\big)-\eta \mb{y}_{k-l}\big\|_2\nonumber\\
\leq&\|\mc{A}_{k-l}-I_{np}\|_{2}\|\widetilde{\mb{x}}^{\textrm{w}}_{k-l}\|_2+\eta \|\mb{y}_{k-l}\|_2
\end{align*}
Summing over~$l$ leads to
{\small \begin{align*}
\sum_{l=0}^{\overline{C}-1}\|\mathbf{x}_{k-(l-1)}-\mathbf{x}_{k-l}\|_2\leq &2\sqrt{n}\|\widetilde{\mb{x}}^{\textrm{w}}_{k-(\overline{C}-1)}\|_2+\eta \|\mb{y}_{k-(\overline{C}-1)}\|_2\\&+\sum_{l=0}^{\overline{C}-2}\Big(2\sqrt{n}\|\widetilde{\mb{x}}^{\textrm{w}}_{k-l}\|_2+\eta\|\mb{y}_{k-l}\|_2\Big).
\end{align*}}From Lemma~\ref{lem:y}, we have~${\forall\, k\geq \overline{C}-1}$:
\begin{align*}
\|\widetilde{\mathbf{s}}^{\textrm{w}}_{k+1}\|_2\leq& m(2\sqrt{n}+\eta Ln)\,(\|\widetilde{\mb{x}}^{\textrm{w}}_{k-(\overline{C}-1)}\|_2+\textstyle \sum_{l=0}^{\overline{C}-2}\|\widetilde{\mb{x}}^{\textrm{w}}_{k-l}\|_2)\\
&+\eta nmL\,\,(\|{\mb{r}}_{k-(\overline{C}-1)}\|_2+\textstyle\sum_{l=0}^{\overline{C}-2}\|{\mb{r}}_{k-l}\|_2)\\
&+(\eta m+{\gamma_{\mc{B}}})\|\widetilde{\mathbf{s}}^{\textrm{w}}_{k-(\overline{C}-1)}\|_2+\eta m\textstyle\sum_{l=0}^{\overline{C}-2}\|\widetilde{\mb{s}}^{\textrm{w}}_{k-l}\|_2,
\end{align*}
and the lemma follows.
\end{IEEEproof}

\subsection{Proof of Lemma~\ref{lem:M_0}}\label{M0_ev}
\begin{IEEEproof} 
Based on the successive application of Schur's determinant identity~\cite{matrix,powell2011calculating}, the characteristic polynomial of~$M^0$ is given by~$(-1)^{\overline{C}}(\lambda-1)(\gamma_{\mc{A}}-\lambda^{\overline{C}})(\gamma_{\mc{B}}-\lambda^{\overline{C}})\lambda^{(\overline{C}-1)}.$
Given that~$\gamma_{\mc{A}},\gamma_{\mc{B}}\in(0,\,1)$, the spectral radius of~$M^0$ is~$1$ and the corresponding eigenvalue,~$\lambda=1$, is simple. We now proceed to find the corresponding right and left eigenvectors,~$\mb{u}$ and~$\mb{w}$, respectively.

By decomposing~$\mb{u}$ and~$\mb{w}$ as follows,
\begin{align*}
\mb{u} = \begin{bmatrix}
{\mb{u}^1}\\
\vdots\\
{\mb{u}^{\overline{C}}}
\end{bmatrix}
,\qquad
\mb{w} = \begin{bmatrix}
{\mb{w}^1}\\
\vdots\\
{\mb{w}^{\overline{C}}}
\end{bmatrix},
\end{align*}
where each~$\mb{u}^i,\mb{w}^i$ is in~$\mbb{R}^3$, from~$M^0\mb{u}=\mb{u}$, we get
\begin{align*}
\gamma_A&[\mb{u}^{\overline{C}}]_1=[\mb{u}^{1}]_1,&{\mb{u}^i=\mb{u}^{i+1},\quad 1\leq i \leq \overline{C}-1},
\end{align*}
resulting in~$[\mb{u}^{i}]_1=0,\,\forall\,i$. Furthermore,
\begin{align*}
\gamma_B&[\mb{u}^{\overline{C}}]_3+2m\sqrt n\sum_{i=1}^{\overline{C}}[\mb{u}^{i}]_1=[\mb{u}^{1}]_{3}.
\end{align*}
Therefore,~$\gamma_B[\mb{u}^{\overline{C}}]_3=[\mb{u}^{1}]_3$ which implies~$[\mb{u}^{i}]_3=0,\, \forall\,i$.
The entries~$[\mb{u}^{i}]_2$ are free variables and we set them equal to~${1,\,\forall\,i}$. Consequently,~${{\mb{u}^i}=\underline{\mb{u}}=\begin{bmatrix}0&1&0\end{bmatrix}^{\top}}$
and~$\mb{u}=\mb{1}_{\overline{C}}\otimes \underline{\mb{u}}$.

Similarly, from~$\mb{w}^{\top}M^0=\mb{w}^{\top}$, we have
\begin{align*}
\begin{matrix}
2m\sqrt[]{n}[\mb{w}^{1}]_3+[\mb{w}^{i+1}]_1=[\mb{w}^{i}]_1,&\qquad i= 1,\ldots,\overline C-1,\\
2m\sqrt[]{n}[\mb{w}^{1}]_3+\gamma_A[\mb{w}^{1}]_1=[\mb{w}^{\overline{C}}]_1.
\end{matrix}
\end{align*}
Summing over all~$i$, we obtain
\begin{align}
2\overline{C}m\sqrt[]{n}[\mb{w}^{1}]_3+\gamma_A[\mb{w}^{1}]_1=[\mb{w}^{1}]_1.\label{eq:qqbangbang}
\end{align}
Furthermore,
\begin{align*}
\begin{matrix}
[\mb{w}^{1}]_2+[\mb{w}^{2}]_2=[\mb{w}^{1}]_2,\\
[\mb{w}^{3}]_2=[\mb{w}^2]_2,\\
\vdots\\
[\mb{w}^{\overline{C}}]_2=[\mb{w}^{\overline{C}-1}]_2,\\
[\mb{w}^{\overline{C}}]_2=0,
\end{matrix}
\end{align*}
resulting in~$[\mb{w}^{i}]_2=0, \forall\, 2\leq i\leq \overline{C}$. Note that~$[\mb{w}^{1}]_2$ is a free variable and we can set it equal to~$1$. Additionally,
\begin{align*}
&\begin{matrix}
[\mb{w}^{2}]_3=[\mb{w}^{1}]_3,\\
[\mb{w}^{3}]_3=[\mb{w}^2]_3,\\
\vdots\\
[\mb{w}^{\overline{C}}]_3=[\mb{w}^{{\overline{C}-1}}]_3,\\
\gamma_B[\mb{w}^{1}]_3=[\mb{w}^{\overline{C}}]_3,
\end{matrix}
\end{align*}
resulting in~$[\mb{w}^{i}]_3=0$, which from~Eq.~\eqref{eq:qqbangbang} implies~${[\mb{w}^{i}]_1=0}$, for all~$i$. Consequently,
\begin{align*}
\mb{w}^{\top}=&\begin{bmatrix}
0 & [\mb{w}^{1}]_2=1&0&0&\cdots&0
\end{bmatrix}.
\end{align*}
\end{IEEEproof}

{\small
\bibliographystyle{IEEEtran}
\bibliography{sample}
}

\end{document}